\documentclass[12pt,a4paper,reqno]{amsart}
\usepackage{amsmath}
\usepackage{amscd}
\usepackage{amssymb}
\usepackage[hypertex]{hyperref}%
\usepackage{latexsym}

\theoremstyle{plain} %% This is the default
\newtheorem{thm}{Theorem}[section]
\newtheorem{cor}[thm]{Corollary}
\newtheorem{lem}[thm]{Lemma}
\newtheorem{prop}[thm]{Proposition}

\newtheorem{defn}[thm]{Definition}
\newtheorem{ex}[thm]{Example}

\newtheorem{rk}[thm]{Remark}

\theoremstyle{remark}

\numberwithin{equation}{section}

\newcommand{\thmref}[1]{Theorem~\ref{#1}}
\newcommand{\propref}[1]{Proposition~\ref{#1}}

\newcommand{\lemref}[1]{Lemma~\ref{#1}}
\newcommand{\corref}[1]{Corollary~\ref{#1}}

\newcommand{\exref}[1]{Example~\ref{#1}}

\newsavebox{\SmallMathBox}

\def\pdo{\Psi{\rm DO}}

\def\Ci{C^\infty}

\def\dd{\partial}

\def\Di{D\kern -.65em /}
\def\Dii{D\kern -.45em /}
\def\di{{\dd}\kern -.55em /}
\def\dii{{\dd}\kern -.40em /}

\def\noi{\noindent}

\def\to{\rightarrow}
\def\too{\longrightarrow}
\def\mto{\mapsto}
\def\mtoo{\longmapsto}

\def\Dd{{\mathcal D}}

\def\Nn{{\mathcal N}}

\def\Rr{{\mathcal R}}

\def\re{{\rm Re}}

\def\={\cong}
\def\>{\supset}
\def\<{\subset}
\def\ii{^{-1}}

\def\12{\frac{1}{2}}
\def\2{\Dd}
\def\3{\Nn}
\def\4{\Rr}
\def\6{\cup}
\def\8{\otimes}
\def\0{^{\circ}}
\def\){\hfill{\ \qed}\enddemo}

\def\a{\alpha}

\def\b{\beta}
\def\R{\mathbb{R}}
\def\C{\mathbb{C}}

\def\e{\varepsilon}

\def\g{\gamma}
\def\G{\Gamma}

\def\th{\theta}

\def\k{\kappa}

\def\la{\lambda}

\def\N{\NN}
\def\o{\infty}
\def\Si{\Sigma}
\def\z{\zeta}
\def\Z{\ZZ}

\def\Cl{\mbox{\rm Cl}}
\def\CS{\mbox{\rm CS}}
\def\Ell{\mbox{\rm Ell}}

\def\Dd{{\mathcal D}}

\def\Nn{{\mathcal N}}

\def\Rr{{\mathcal R}}

\def\index{\mbox{\rm index\,}}

\def\Si{S\kern -.65em /}

\def\supp{\mbox{\rm supp}}

\def\tr{\mbox{\rm tr\,}}
\def\Tr{\mbox{\rm Tr\,}}

\def\TRx{\textmd{\Small {\rm TR}}_x}

\def\res{\textmd {\rm res}}
\def\resx{\textmd {\rm res}_x}
\def\resxo{\textmd {\rm res}_{x,0}}

\def\dbar{d{\hskip-1pt\bar{}}\hskip1pt}

\def\Cf{\mathbb{C}}

\def\Nf{\mathbb{N}}

\def\Rf{\mathbb{R}}

\def\Zf{\mathbb{Z}}

\def\ord{\mbox{\rm ord}}

\def \C{{\! \rm \ I \!\!\!C}}
\def \R {{\! \rm \ I \!R}}
\def \N {{\! \rm \ I \!N}}

\def \Z {{\! \rm Z\! \!Z}}

\def\cutoffint{-\hskip -12pt\int}

\def\otherterm#1{{\it#1}}
\def \e {{\epsilon}}

\def \s {{\sigma}}

\def \End {{\rm End}}

\def \ord {{\rm ord}}
\def \Ci {{C^\infty}}
\def \l {{\lambda}}
\def \tr {{\rm tr}}
\def\pdo{\psi{\rm do}}

\def\TR{{\rm TR}}

\def\res{{\rm res}}

\title[running header]{ cover page title of paper}
\begin{document}

\title[]
{\bf A Laurent expansion for regularized integrals of holomorphic
symbols \vskip 7mm  {\rm {\small Sylvie PAYCHA {\tiny and\,} Simon
SCOTT}}\vskip 5mm}

\begin{abstract}  For a  holomorphic family  of classical
pseudodifferential operators on a closed manifold we give exact
formulae for all coefficients in the Laurent expansion of its
Kontsevich-Vishik canonical trace. This generalizes to all higher
order terms a known result identifying the residue trace with a
pole of the canonical trace.
\end{abstract}

\maketitle

\section*{Introduction}
Let $M$ be a compact boundaryless Riemannian manifold of dimension
$n$ and $E$ a smooth vector bundle based on $M$. For a classical
pseudodifferential operator ($\pdo$) $A$ with non-integer order
acting on smooth sections of $E$ one can define following
Kontsevich and Vishik \cite{KV} and Lesch \cite{Le} the canonical
trace of $A$
$${\rm TR}(A):=\int_M dx\ \TRx(A)\ ,
\quad \TRx(A):= \cutoffint_{T_x^*M} {\rm tr}_x(\sigma_A(x, \xi))\
\dbar\xi,$$ in terms of a local classical symbol $\sigma_A$ and a
finite-part integral $\cutoffint_{T_x^*M}$ over the cotangent
space $T_x^*M$ at $x\in M$. Here, $\dbar\xi = (2\pi)^{-n}\,d\xi$
with $d\xi$ Lebesgue measure on $T_x^*M\cong\Rf^n$, while ${\rm
tr}_x$ denotes the fibrewise trace. Since the work of Seeley
\cite{Se1} and later of Guillemin \cite{Gu}, Wodzicki \cite{Wo}
and then Kontsevich and Vishik \cite{KV}, it has been known that
given a holomorphic family $z\mapsto A(z)$  of classical $\pdo$s
parametrized by  a domain $W\subset \C$, with holomorphic order
$\alpha:W\to \Cf$ such that $\alpha^\prime $ does not vanish on
$$P:= \alpha^{-1}\left(\Z\cap [-n, +\infty[\right) ,$$ then the
map $z\mapsto {\rm TR}(A(z))$ is a meromorphic function with no
more than simple poles located in $ P$. The complex residue at
$z_0\in P$ is given by  a local expression \cite{Wo}, \cite{Gu},
\cite{KV}
\begin{equation}\label{e:resTR}
{\rm Res}_{z=z_0}\,{\rm TR}(A(z))
=-\frac{1}{\alpha^\prime(z_0)}\,{\rm res}(A(z_0)),
\end{equation}
where for a classical pseudodifferential operator $B$ with symbol
$\sigma_B$
$${\rm res}(B):=  \int_M dx\ { \rm res}_x(B), \quad {\rm
res}_x(B):=\int_{S_x^*M} {\rm tr}_x\left(\,(\sigma_B)_{-n}(x,
\xi)\right)\, \dbar_S\xi$$ is the residue trace of $B$. Here,
$\dbar_S\xi = (2\pi)^{-n}\,d_S(\xi)$ with $d_S(\xi)$ the sphere
measure on $S_x^*M = \{|\xi|=1 \ | \ \xi \in T_x^*M\}$, while the
subscript refers to the positively homogeneous component of the
symbol of order $-n$.

\vskip 2mm

In this paper, extending the identification \eqref{e:resTR}, we
provide a complete solution to the problem of giving exact
formulae for all coefficients in the Laurent expansion of ${\rm
TR}(A(z))$ around each pole in terms of locally defined canonical
trace and residue trace densities.

\vskip 1mm

 For a meromorphic function $G$, define its finite-part
${\rm fp}_{z=z_0}G(z)$ at $z_0$ to be the constant term  in the
Laurent expansion of $G(z)$ around $z_0$. Let $A^{(r)}(z) =
\dd^r_z A(z)$ be the derivative $\pdo$ with symbol
$\s_{A^{(r)}(z)} := \dd_z^r\s_{A(z)}$.

\vskip 5mm

\noi {\bf Theorem}  \ {\it Let $z\mto A(z)$ be a holomorphic
family of classical $\pdo$s of order $\a(z) = qz + b$. If $z_0\in
P$ and $q\neq 0$, then ${\rm TR}(A(z))$ has Laurent expansion for
$z$ near $z_0$
\begin{equation}\label{eq:intro1}
{\rm TR}(A(z))\ = \ -\,\frac{{\rm
res}(A(z_0))}{q}\,\frac{1}{(z-z_0)} \ + \ \sum_{k=0}^{\o}{\rm
fp}_{z=z_0} {\rm TR}(A^{(k)}(z))\frac{(z-z_0)^k}{k!}.
\end{equation}
Furthermore,
\begin{equation}\label{eq:density}
\left(\TRx \left(A^{(k)}(z_0)\right)  \ - \ \frac{1}{q\,
(k+1)}\,\resxo\left(A^{(k+1)}(z_0)\right)\right) \ dx\,
\end{equation}
\noi defines a global density on $M$ and
\begin{equation}\label{eq:intro2}
{\rm fp}_{z=z_0} {\rm TR}(A^{(k)}(z))  \ = \  \int_M dx\,
\left(\TRx \left(A^{(k)}(z_0)\right)  \ - \ \frac{1}{q\,
(k+1)}\,\resxo\left(A^{(k+1)}(z_0)\right)\right) \ .
\end{equation}
At a point  $z_0\notin P$ the function ${\rm TR}(A(z))$ is
holomorphic near $z_0$ and the Laurent expansion \eqref{eq:intro1}
reduces to the Taylor series
\begin{equation*}
{\rm TR}(A(z))\ = \ {\rm TR}(A(z_0)) \ + \ \sum_{k=1}^{\o}{\rm
TR}(A^{(k)}(z_0))\,\frac{(z-z_0)^k}{k!}.
\end{equation*}}

\vskip 2mm \noi It is to be emphasized here that $A^{(r)}(z)$
cannot be a classical $\pdo$ for $r> 0$, but in local coordinates
is represented for $|\xi|>0$ by a log-polyhomogeneous symbol of
the form
$$\s_{A^{(r)}(z)}(x,\xi) \ \sim \
\sum_{j\geq 0}\sum_{l=0}^r \s(A^{(r)}(z))_{\alpha(z)-j,\,
l}(x,\xi)\,\log^l |\xi|$$ with $\s(A^{(r)})_{\alpha(z)-j,
\,l}(x,\xi)$ positively homogeneous in $\xi$ of degree $\a(z)-j$.
It follows that individually the terms in \eqref{eq:density}
\begin{equation}\label{eq:TRx}
\TRx \left(A^{(k)}(z_0)\right) \, dx :=
\cutoffint_{T_x^*M} {\rm tr}_x(\s_{A^{(k)}(z_0)}(x,\xi))\
 \dbar\xi \, dx
\end{equation}
and
\begin{equation}\label{eq:resx}
\resxo\left(A^{(k+1)}(z_0)\right) \, dx := \int_{S_x^*M} {\rm
tr}_x\left(\,\s(A^{(k+1)}(z_0))_{-n,0}(x,\xi)\right)\,
\dbar_S\xi\, dx
\end{equation}
do not in general  determine globally defined densities on the
manifold $M$ when $r>0$, rather it is then only the sum of terms
\eqref{eq:density} which integrates to a global invariant of $M$.
(In particular, it is important to distinguish \eqref{eq:resx}
from the higher residue trace density of \cite{Le}, see
Remark(\ref{higherresidue}) here.) When $\a(z) = qz + b$ is {\em
not} integer valued it is known that $\TRx(A^{(k)}(z))\,dx$ does
then define a global density on $M$; in this case,
$\resxo(A^{(k+1)}(z))$ is identically zero and \eqref{eq:intro2}
reduces to  the canonical trace $\TR(A^{(k)}(z)) = \int_M dx\,
\TRx(A^{(k)}(z))$ on non-integer order $\pdo$s with
log-polyhomogeneous symbol \cite{Le}.

\vskip 1mm

These results hold more generally when $\a(z)$ is an arbitrary
holomorphic function with $\a^\prime(z_0)\neq 0$ at $z_0\in P$.
Then the local residue term in \eqref{eq:intro2} is replaced by
the local residue of an explicitly computable polynomial in the
symbols of the operators $A^{(k+1)}(z_0), \ldots, A(z_0)$. A
general formula is given in \thmref{thm:KV2}, here we state the
formula just for the constant term in the Laurent expansion of
${\rm TR}(A(z))$: one has

\begin{equation}\label{eq:intro2a}
{\rm fp}_{z=z_0} {\rm TR}(A(z))  \ = \  \int_M dx\, \left(\TRx
\left(A(z_0)\right)  \ - \
\frac{1}{\alpha^\prime(z_0)}\,\resxo\left(A^{'}(z_0)\right)\right)
\end{equation}
$$\hskip 10mm + \ \
\frac{\a^{\prime\prime}(z_0)}{2\,\a^\prime(z_0)^2}\
\res\,\left(A(z_0) \right).$$

\vskip 3mm

\noi Thus compared to \eqref{eq:intro2}, the constant term
\eqref{eq:intro2a} in the expansion acquires an additional residue
trace term. Moreover, the identification implies that

\begin{equation}\label{eq:density0}
\left(\TRx \left(A(z_0)\right)  \ - \
\frac{1}{\alpha^\prime(z_0)}\,\resxo\left(A^\prime(z_0)\right)\right)
\ dx\,
\end{equation}

\vskip 3mm

\noi defines a global density on $M$ independently of the order
$\a(z_0)$ of $A(z_0)$ (for $z_0\notin P$ the locally defined
residue term vanishes and \eqref{eq:density0} reduces to the usual
canonical trace density). Though this follows from general
properties of holomorphic families of canonical traces, we
additionally give an elementary direct proof in Appendix A.

\vskip 4mm

Applied to $\pdo$ {\em zeta-functions} this yields formulae for a
number of widely studied spectral geometric invariants. For $Q$ an
elliptic classical pseudodifferential operator of order $q>0$ and
with spectral cut $\theta$, its complex powers $Q_\theta^{-z}$ are
well defined \cite{Se1} and to a classical pseudodifferential
operator $A$ of order $\a\in\R$ one can associate the holomorphic
family $A(z)=A\, Q_\theta^{-z}$ with order function $\a(z) = \a -q
z$. The generalized zeta-function $$z\mapsto \zeta_\theta(A, Q,
z):= {\rm TR}(A\, Q_\theta^{-z})$$ is meromorphic on $\Cf$ with at
most simple poles in $ P:=\{\frac{\alpha-j}{q}\ | \ j\in [-n,
\infty) \ \cap\ \Z\}$. It has been shown by Grubb and Seeley
\cite{GS,Gr1} that $\G(s)\, \zeta_\theta(A,Q, s)$ has pole
structure
\begin{equation}\label{e:zetapoles}
\G(s)\, \zeta_\theta(A,Q, s) \sim \sum_{j\geq -n} \frac{c_j}{s +
\frac{j - \a}{q}} -  \frac{\Tr(A \, \Pi_Q)}{s} + \sum_{l\geq 0}
\left(\frac{c_l^{'}}{(s+l)^2} + \frac{c_l^{''}}{(s+l)}\right),
\end{equation}
where the coefficients $c_j$ and $c_l^{'}$ are locally determined,
by finitely many homogeneous components of the local symbol, while
the $c_l^{''}$ are globally determined. In particular, whenever
$$\frac{j-\a}{q} :=l \in [0,\infty) \ \cap\ \Z$$ it is shown that
the sum of terms
\begin{equation}\label{e:globalsumpoles}
c_l^{''} + c_{\a + l q}
\end{equation}
is defined invariantly on the manifold $M$, while individually the
coefficients $c_l^{''}$ and $c_{\a + l q}$ (which contain
contributions from the terms \eqref{eq:TRx} and \eqref{eq:resx}
respectively) depend on the symbol structure in each local
trivialization. Here, in \thmref{thm:KV3}, we compute the Laurent
expansion around each of the poles of the meromorphically
continued Schwartz kernel $K_{A\,Q_{\theta}^{-z}}(x,x)\vert^{{\rm
mer}}:= \cutoffint_{T_x^*M} \sigma_{A\, Q_\theta^{-z}}(x, \xi) \,
\dbar\xi$ giving the following exact formula for
\eqref{e:globalsumpoles}. One has

\begin{equation}\label{e:globalsumpolesformula}
c_{l}^{''} + c_{\a + l q}   = \frac{(-1)^l}{l!}\int_M dx \left(
\TRx(A\,Q^l) \ - \ \frac{1}{q} \, \resxo(A\,Q^l\,\log_\theta
Q)\right).
\end{equation}

\vskip 2mm

\noi The remaining coefficients in \eqref{e:zetapoles} occur as
residue traces of the form \eqref{e:resTR}. By a well known
equivalence, see for example \cite{GS,Gr1}, when $Q$ is a
Laplace-type operator these formulae acquire a geometric character
as coefficients in the asymptotic heat trace expansion
 $\Tr(A\,e^{-tQ}) \sim \sum_{j\geq -n} c_j
\,t^\frac{j - \a}{q} + \sum_{l\geq 0} (-\,c_l^{'} \,\log t +
c_l^{''})\, t^l$ as $t\to 0+$.

\vskip 1mm

From \eqref{e:resTR} (\cite{Wo,Gu,KV}) $\zeta_\theta(A,Q, z)$ has
a simple pole at $z=0$ with residue $-\frac{1}{q}\,{\rm res}(A)$,
which vanishes if $\a\notin\Z$. The coefficients of the full
Laurent expansion of $\zeta_\theta(A,Q, s)$ around $z=0$ are given
by the following formulae (Theorem \ref{thm:zetamer}).

\vskip 4mm

\noi {\bf Theorem} \ {\it For $k\in \Nf$, let
$\zeta_\theta^{(k)}(A,Q, 0)$ denote the coefficient of $z^k/k!$ in
the Laurent expansion of $\zeta_\theta(A,Q, z)$ around $z=0$. Then
\begin{eqnarray}
\zeta_\theta^{(k)}(A,Q, 0) \ & = &\ (-1)^k\, \int_M dx\,
\left(\TRx (A\, \log_\theta^k Q) - \frac{1}{q\,
(k+1)}\,\resxo(A\,\log_\theta^{k+1} Q )\right) \
\nonumber \\[2mm] && \hskip 20mm + \ \  (-1)^{k+1}\, {\rm tr}\left(A\,
\, \log_\theta^k Q\, \Pi_Q\right) \ , \label{eq:intro3b}
\end{eqnarray}

\vskip 2mm \noi where $\Pi_Q$ is a smoothing operator projector
onto the generalized kernel of $Q$. Specifically, for a classical
$\pdo$ $A$ of arbitrary order
\begin{equation}\label{eq:density3a}
\left(\TRx (A) - \frac{1}{q}\ \resxo(A\,\log_\theta Q )\right)\,dx
\end{equation}
is a globally defined density on $M$ and, setting
$\zeta_\theta(A,Q, 0):= \zeta_\theta^{(0)}(A,Q, 0)$,  the constant
term in the expansion around $z=0$ is
\begin{equation}\label{eq:intro3a} \zeta_\theta(A,Q,
0)\ = \ \int_M dx\, \left(\TRx (A) - \frac{1}{q}\
\resxo(A\,\log_\theta Q )\right)  - {\rm
tr}\left(A\, \Pi_Q\right) \ .\\[1mm]
\end{equation}
\noi When $A$ is a differential operator $\zeta_\theta(A, Q, 0)=
\lim_{z\to 0} \,\zeta_\theta(A, Q, z)$ and equation
\eqref{eq:intro3a} becomes
\begin{equation}\label{eq:intro4}
\zeta_\theta(A, Q, 0) \ = \ -\frac{1}{q}\ {\rm res}(A\,
\log_\theta Q)-{\rm tr}\left(A\, \Pi_Q\right) \ .
\end{equation}
\noi When $Q$ is a differential operator and $m$ a non-negative
integer, setting $\zeta_\theta(Q, -m):= {\rm
fp}_{z=-m}\zeta_\theta(I,Q,z)$, one has
\begin{equation}\label{eq:intro4a}
\zeta_\theta(Q, -m) \ = \ -\frac{1}{q}\ {\rm res}(Q^m\,
\log_\theta Q)-{\rm tr}\left(Q^m\, \Pi_Q\right) \ .
\end{equation}
\vskip 1mm

 If $A$ is a $\pdo$ of non-integer order $\a\notin\Z$ then
$0\notin P$ and from \cite{Le} the canonical trace of $A\,
\log_\theta^k Q$ is defined. Then \eqref{eq:intro3b} reduces to
\begin{equation}\label{eq:intro3b2}
\zeta_\theta^{(k)}(A,Q, 0) =  (-1)^k\,  \TR(A\, \log_\theta^k Q) -
(-1)^k\, {\rm tr}\left(A\, \, \log_\theta^k Q\, \Pi_Q\right)
\end{equation}
and, in particular, in this case
\begin{equation}\label{eq:intro3b3}
\zeta_\theta (A,Q, 0) =  \TR (A)  - \, {\rm tr}\left(A\,
\Pi_Q\right).
\end{equation}
}

\vskip 2mm

\noi  Notice, that in \eqref{eq:intro4} the term ${\rm res}(A\,
\log_\theta Q) =\zeta_\theta(A, Q, 0)+{\rm tr}\left( A\,
\Pi_Q\right)$ is locally determined, depending on only finitely
many of the homogeneous terms in the local symbols of $A$ and $Q$
(\cite{GS} Thm 2.7, see also \cite{Sc} Prop 1.5). In the case
$A=I$ the identity \eqref{eq:intro4} was shown for
pseudodifferential $Q$ in \cite{Sc} and \cite{Gr2}, and in the
particular case where $Q$ is an invertible positive differential
operator \eqref{eq:intro4a} can be inferred from \cite{Lo}. The
identity \eqref{eq:intro3b3} is known from \cite{Gr1}(Rem. (1.6)).
A resolvent proof of \eqref{eq:intro3a} has been given recently in
\cite{Gr3}. \vskip 2mm

If, on the other hand, one considers, for example, $A(z) = A\,
Q_\theta^{-\frac{z}{1+\mu \, z}}$, then  the corresponding `zeta
function' $\TR(A\, Q_\theta^{-\frac{z}{1+\mu \, z}})$ has simple
and real poles in $\Cf\backslash\{-1/\mu\}$ and by
\eqref{eq:intro2a} the constant term at $z=0$  has, compared to
\eqref{eq:intro3a}, an extra term
\begin{equation*}\label{eq:genz}
{\rm fp}_{z=0} \TR(A\, Q_\theta^{-\frac{z}{1+\mu \, z}}) \ = \
\int_M dx\, \left(\TRx (A)\ -\  \frac{1}{q}\ \resxo(A\,\log_\theta
Q )\right) \ -\  {\rm tr}\left(A\, \Pi_Q\right)
\end{equation*}
$$+\ \
\frac{\mu}{q}\ \res\,\left(A\right).$$ \vskip 1mm \noi The
appearance here of $\frac{\mu}{q}\ \res\,\left(A\right)$
corresponds to additional terms that occur as a result of a
rescaling of the cut-off parameter when expectation values are
computed from Feynman diagrams using a momentum cut-off procedure,
see \cite{Gro}. See also Remark \ref{rk}.

\vskip 1mm

One view point to adopt on \eqref{eq:intro2a} is that it provides
a defect formula for regularized traces and indeed most well known
trace defect formulas \cite{MN,O1,CDMP,Gr2} are an easy
consequence of it. On the other hand, new more precise formulae
also follow. In particular, though $\TR$ is not in general defined
on the bracket $[A,B]$ when the bracket is of integer order, we
find (Theorem \ref{e:TRacedefect}) that in this case the following
exact global formula holds
\begin{equation}\label{e:TRbracketnoninteger}
\int_M dx\, \left(\,\TRx\left([A,B]\right) -
\frac{1}{q}\,\resxo\left([A, B\log_\theta Q]\right)\,\right) \ = \
0
\end{equation}
independently of the choice of  $Q$; when $[A,B]$ is not of
integer order \eqref{e:TRbracketnoninteger} reduces to the usual
trace property of the canonical trace $\TR([A,B])=0$, see Section
2.

\vskip 2mm

Looking at the next term up in the Laurent expansion of
$\TR(Q^{-z})$ at zero, equation (\ref{eq:intro3b}) provides an
explicit formula  for the $\zeta$-determinant
$${\rm det}_{\zeta, \theta} Q \ = \exp(-\,\zeta_\theta^\prime(Q, 0)),$$
where $\zeta_\theta^\prime(Q, 0) = \dd_z\zeta_\theta(Q,
z))_{\vert_{z=0}}$, of an invertible elliptic classical
pseudodifferential operator $Q$ of positive order $q$ and with
spectral cut $\theta$. The zeta determinant is a complicated
non-local invariant which has been studied in diverse mathematical
contexts. From \eqref{eq:intro3b} one finds (\thmref{t:det}):

\vskip 2mm

\noi {\bf Theorem} \
\begin{equation}\label{eq:logdet}
\log {\rm det}_{\zeta, \theta} (Q) = \int_M dx\left( \TRx
\left(\log_\theta Q\right)-\frac{1}{ 2q}\, \resxo
\left(\log^2_\theta Q\right)\right).
\end{equation}

\vskip 2mm \noi A slightly modified formula holds for
non-invertible $Q$. Notice, here, that $\TR$ of $\log_{\th}Q$ does
not generally exist; if $\resxo(\log_{\th}^2 Q)=0$ pointwise it is
defined, and then $\log {\rm det}_{\zeta, \theta} (Q) =
\TR(\log_{\th}Q)$, which holds for example for odd-class operators
of even order, such as differential operators of even order on
odd-dimensional manifolds \cite{KV,O2}. Equation \eqref{eq:logdet}
leads to explicit formulae for the multiplicative anomaly.

\vskip 5mm

\section{Finite-part integrals (and canonical traces) of holomorphic
families of classical symbols (and pseudodifferential operators)}

\subsection{Classical and log-polyhomogeneous symbols} We briefly
recall some  notions concerning symbols and pseudodifferential
operators and fix the corresponding notations. Classical
references for the polyhomogeneous symbol calculus are e.g.
\cite{Gi}, \cite{GS}, \cite{Ho}, \cite{Se2}, \cite{Sh}, and for
the extension to log-polyhomogeneous symbols \cite{Le}. $E$
denotes a smooth hermitian vector bundle based on some closed
Riemannian manifold $M$. The space $\Ci(M, E)$ of smooth sections
of $E$  is endowed with  the inner product $\langle \psi,
\phi\rangle:= \int_M d \mu(x) \langle \psi(x), \phi(x)\rangle_x$
induced by the hermitian structure $\langle \cdot, \cdot\rangle_x$
on the fibre over $x\in M$  and the Riemannian measure $\mu$ on
$M$. $H^s(M,
E)$ denotes the  $H^s$-Sobolev closure of the space $\Ci(M, E)$.\\

\noi Given an open subset $U$ of $\R^n$ and an auxiliary
(finite-dimensional) normed vector space $V$, the set of symbols
${\rm S}^r(U,V)$ on $U$ of order $r\in\Rf$  consists of those
functions $\s(x,\xi)$ in $\Ci(T^*U,{\rm End}( V))$ such that
$\dd_x^{\mu}\dd_{\xi}^{\nu}\s(x,\xi)$ is $O((1+|\xi|)^{r -
|\nu|})$ for all multi-indices $\mu, \nu$, uniformly in $\xi$,
and, on compact subsets of $U$, uniformly in $x$. We set ${\rm
S}(U, V):= \bigcup_{r\in \R} {\rm S}^r(U,V)$ and ${\rm
S}^{-\infty}(U,V):= \bigcap_{r\in \R} {\rm S}^r(U,V)$. A {\it
classical} (1-step polyhomogeneous) symbol of order $\a\in\Cf$
means a  function $\s(x,\xi)$ in $\Ci(T^*U,{\rm End}( V))$  such
that for each $N\in \N$ and each integer $0\leq j \leq N$ there
exists $\sigma_{\alpha-j}\in \Ci(T^*U,{\rm End}( V))$ which is
homogeneous in $\xi$ of degree $\a-j$ for $|\xi|\geq 1$, so
$\sigma_{\alpha-j}(x,t\xi) = t^{\a-j}\, \sigma_{\alpha-j}(x,\xi)$
for $t\geq 1, |\xi|\geq 1$, and a symbol $\sigma_{(N)}\in {\rm
S}^{\,\re(\alpha)-N-1}(U,V)$ such that
\begin{equation}\label{eq:classymb}
\sigma (x, \xi)= \sum_{j=0}^N  \sigma_{\alpha-j}(x,\xi)+
\sigma_{(N)}(x, \xi)\quad \forall  \ (x, \xi)\in T^*U.
\end{equation}
We then write $\sigma(x, \xi)\sim \sum_{j=0}^\infty
\sigma_{\alpha-j}(x,\xi).$ Let $\CS(U,V)$ denote the class of
classical symbols on $U$ with values in $V$ and let $\CS^\alpha
(U,V)$ denote the subset of classical symbols of order $\alpha$.
When $V=\C$, we write ${\rm S}^r(U)$, $\CS^\alpha(U)$, and so
forth; for brevity we may omit the $V$ in the statement of some
results. A $\pdo$ which for a given atlas on $M$ has a classical
symbol in the local coordinates defined by each chart is called
{\it classical}, this is independent of the choice of atlas. Let
$\Cl(M, E)$ denote the algebra of classical $\pdo$s acting on
$\Ci(M, E)$ and let $\Ell(M, E)$ be the subalgebra of elliptic
operators. For any $\alpha \in \C$ let $\Cl^\alpha(M, E)$, resp.
$\Ell^\alpha(M, E)$, denote the subset of operators in $\Cl(M,E)$,
resp.  $\Ell (M, E)$, of order $\alpha$. With $\R_+ = (0,\o)$, set
$\Ell_{\ord
>0} (M,E):= \bigcup_{r\in\R_+}\Ell^r(M, E)$.

\vskip 2mm

To deal with derivatives of complex powers of classical $\pdo$s
one considers the larger class of $\pdo$s with log-polyhomogeneous
symbols. Given an open subset $U\subset M$, a non-negative integer
$k$ and a complex  number $\alpha$, a symbol $\sigma$ lies in
$\CS^{\alpha, k}(U,V)$ and is said to have order $\alpha$ and log
degree $k$ if
\begin{equation}\label{eq:logpolysymb} \sigma (x, \xi)=
\sum_{j=0}^N \sigma_{\alpha-j}(x,\xi)+ \sigma_{(N)}(x, \xi)\quad
\forall  \ (x, \xi)\in T^*U
\end{equation}
where $\sigma_{(N)}\in {\rm S}^{\re(\alpha)-N-1+\e}(U,V)$ for any
$\e>0$,  and  $$ \sigma_{\alpha-j}(\xi) = \sum_{l=0}^k
\sigma_{\alpha-j,\, l} ( x,\xi)  \log^l [\xi] \quad \forall \
\xi\in T_x^*U$$ with $\sigma_{\alpha-j,\, l}$ homogeneous in $\xi$
of degree $\alpha-j$ for $|\xi|\geq 1$, and (in the notation of
\cite{Gr1}) $[\xi]$ a strictly positive $C^{\o}$ function in $\xi$
with $[\xi]=|\xi|$ for $|\xi|\geq 1$. As before, in this case we
write
\begin{equation}\label{asymplogphgs}
\sigma(x, \xi)\sim \sum_{j=0}^\infty \sigma_{\alpha-j}(x,\xi) =
\sum_{j=0}^\infty \sum_{l=0}^k
\sigma_{\alpha-j,\,l}(x,\xi)\,\log^l[\xi].
\end{equation} Then $\CS^{*, *}(U,V):=
\bigcup_{k=0}^\infty   \CS^{*, k}(U,V)$, where $\CS^{*, k}(U,V)=
\bigcup_{\alpha\in \C} \CS^{\alpha, k}(U,V),$ defines the class
filtered by $k$ of log-polyhomogeneous symbols on $U$. In
particular,  $\CS(U,V)$ coincides with $\CS^{*, 0}(U,V)$.

\vskip 1mm

Given a non-negative integer $k$, let $\Cl^{\alpha, k}(M, E)$
denote the space of  pseudodifferential operators on $\Ci(M, E)$
which in any local trivialization $E_{|U} \cong U \times V$ have
symbol in $\CS^{\alpha, k}(U,V)$. Set $\Cl^{*, k}(M, E):=
\bigcup_{\a \in\C}\Cl^{\a, k}(M, E).$
%\end{itemize}
%\end{rk}

\vskip 1mm

The following subclasses of symbols and $\pdo$s will be of
importance in what follows.
\begin{defn}
A log-polyhomogeneous symbol \eqref{asymplogphgs} with integer
order $\a\in\Z$ is said to be {\it even-even} (or, more fully, to
have {\it even-even alternating parity}) if for each $j\geq 0$
\begin{equation}\label{regularparity}
\sigma_{\alpha-j,\,l}(x,-\,\xi) = (-1)^{\alpha-j}\
\sigma_{\alpha-j,\,l}(x,\xi) \ \ \ {\rm for} \ |\xi|\geq 1,
\end{equation}
and the same holds for all derivatives in $x$ and $\xi$. It is
said to be {\it even-odd} (or, more fully, to have {\it even-odd
alternating parity}) if for each $j\geq 0$
\begin{equation}\label{singularparity}
\sigma_{\alpha-j,\,l}(x,-\,\xi) = (-1)^{\alpha-j-1}\
\sigma_{\alpha-j,\,l}(x,\xi) \ \ \ {\rm for} \ |\xi|\geq 1,
\end{equation}
and the same holds for all derivatives in $x$ and $\xi$. A $\pdo$
$A\in \Cl^{\alpha, k}(M, E)$ will be said to be {\it even-even}
(resp. {\it even-odd}) if in each local trivialization any local
symbol $\s_{A}(x,\xi)\in\CS^{\a, k}(U,V)$ representing $A$ (modulo
smoothing operators) has  even-even (resp. even-odd) parity.
\end{defn}
\noi Thus, an even-even symbol with even integer degree  is even
in $\xi$, while an even-odd symbol with even integer degree is odd
in $\xi$; a similar statement holds if the symbol has odd integer
degree.

\begin{rk} The terminology in Definition (1.1) follows
\cite{Gr4}. Kontsevich-Vishik \cite{KV} studied even-even
classical $\pdo$s on odd-dimensional manifolds, calling them {\it
odd-class} operators. Odd-class operators (or symbols) form an
algebra and include differential operators and their parametrices.
The class of operators with even-odd parity symbols on
even-dimensional manifolds, which includes the modulus operator
$|A| = (A^2)^{1/2}$ for $A$ a first-order elliptic self-adjoint
differential operator, was introduced and studied by Grubb
\cite{Gr1}; this class admits similar properties with respect to
traces on $\pdo$s as the odd-class operators, though they do not
form an algebra. In \cite{O2} Okikiolu uses the terminology
`regular parity' and `singular parity' for \eqref{regularparity}
and \eqref{singularparity}.
\end{rk}

\subsection{Finite part integrals of symbols and the canonical trace}

In order to make sense of $\int_{T_x^*M} \sigma (x,\xi)\,
\dbar\xi$ when $\sigma\in \CS^{\alpha, *}(U,V)$ is a
log-polyhomogeneous symbol (the integral diverges a priori if
$\re(\alpha)\geq -n$) on a local subset  $U\subset \R^n$, one can
extract a ``finite part" when $R\to \infty$ from the integral
$\int_{B_x^*(0,R)} \sigma (x,\xi)\, \dbar\xi$   where
$B_x^*(0,R)\subset T_x^*U$ denotes the ball centered at $0$ with
radius $R$   for a given point
$x\in U$.\\

First, though, we introduce the local residue density on
log-polyhomogeneous symbols, which acts as an obstruction to the
finite part integral of a classical symbol defining a global
density on $M$ and measures the anomalous contribution to the
Laurent coefficients at the poles of the finite part integral when
evaluated on holomorphic families of symbols.

\begin{defn}
Given an open subset $U\subset \R^n$, the local Guillemin-Wodzicki
residue is defined for $\sigma\in \CS^{\a}(U,V)$ by
$$
\resx(\sigma)=\int_{S_x^*U}\tr_x\left(\sigma_{-n}(x, \xi)\right)\
\dbar_S\xi$$ and extends to a map $\resxo:\CS^{\a,k}(U,V)\to \Cf$
by the same formula
\begin{eqnarray*}
\resxo(\sigma)&=&\int_{S_x^*U}\tr_x\left(\sigma_{-n}(x,
\xi)\right)\
\dbar_S\xi \\
&=&\sum_{l=0}^k \int_{S_x^*U}\tr_x\left(\sigma_{-n, l} ( x,\xi)
\right) \log^l \vert \xi \vert \ \dbar_S\xi \\
&=&\int_{S_x^*U}\tr_x\left(\sigma_{-n,0}(x, \xi)\right)\
\dbar_S\xi.
\end{eqnarray*}
\end{defn}

\noi When $k>0$ the extra subscript is included in the notation
$\resxo(\sigma)$ as a reminder that it is the residue of the log
degree zero component of the symbol that is being computed. The
distinction is made because when $k>0$ the local densities
$\resxo(\sigma)\,dx$ do not in general define a global density on
$M$, due to cascading derivatives of powers of logs when changing
local coordinates. When $k=0$,  Guillemin \cite{Gu} and Wodzicki
\cite{Wo} showed the following remarkable properties.
\begin{prop}\label{p:resdensity}
Let $A\in \Cl^{\a}(M, E)$ be a classical $\pdo$ represented in a
local coordinate chart $U$ by $\s\in \CS^\a(U,V)$. Then ${\rm
res}_x(\s)\,dx$ determines a global density on $M$, that is, an
element of \, $\Ci(M,|\Omega|)$, which defines the projectively
unique trace on $\Cl^{*,0}(M, E)$.
\end{prop}
Proofs may be found in {\it loc.cit.}, and in Section 2 here. The
first property means that ${\rm res}_x(\s)\,dx$ can be integrated
over $M$. The resulting number
\begin{equation}\label{e:classicalres}
{\rm res}(A):=  \int_M {\rm res}_x(\s)\,dx = \int_M
dx\int_{S_x^*M}\,{\rm tr}_x\left(\sigma_{-n}(x, \xi)\right)
\,\dbar_S\xi
\end{equation}
is known as the residue trace of $A$.  The terminology refers to
the trace property in \propref{p:resdensity} that if the manifold
$M$ is connected and has dimension larger than $1$, then up to a
scalar multiple \eqref{e:classicalres} defines on $\Cl^{*}(M, E)$
the unique linear functional vanishing on commutators
$${\rm res}([A,B]) = 0, \hskip 15mm  A,B\in \Cl^{*}(M, E). $$
Notice, from its definition, that the residue trace also vanishes
on operators of order $<-n$ and on non-integer order operators.

\vskip 1mm

\begin{rk}\label{higherresidue}
The residue trace was  extended by Lesch \cite{Le} to $A\in
\Cl^{\a, k}(M, E)$  with $k>0$ by defining ${\rm res}_k(A):=
(k+1)!\int_M dx\, \int_{S_x^*M}\,{\rm tr}_x\left(\sigma_{-n,k}(x,
\xi)\right) \,\dbar_S\xi.$ For an operator with
log-polyhomogeneous symbol of log degree $k>0$ the form
$\sigma_{-n,k}(x, \xi)\,dx$ defines a global density on $M$, a
property which is not generally true for the lower log degree
densities $\sigma_{-n,0}(x, \xi)\,dx, \ldots, \sigma_{-n,k-1}(x,
\xi)\,dx$ which depend on the symbol structure in each local
coordinate chart. We emphasize that the higher residue is not
being used in the Laurent expansions we compute here, rather the
relevant object is the locally defined form $\sigma_{-n,0}(x,
\xi)\,dx$ which for suitable $A$ defines one component of a
specific local density which does determine an element of
$\Ci(M,\End(E)\otimes|\Omega|).$
\end{rk}

\vskip 2mm

It was, on the other hand, observed by Kontsevich and Vishik
\cite{KV} that the usual $L^2$-trace on $\pdo$s of real order
$<-n$ extends to a functional on the space
$\Cl^{\C\backslash\Z}(M,E)$ of $\pdo$s of  non-integer order and
vanishes on commutators of non-integer order. Lesch \cite{Le}
subsequently showed that the resulting canonical trace can be
further extended to
$$\Cl^{{\C\backslash\Z}, *}(M, E):=
\bigcup_{\a\, \in \, \C\backslash\Z} \Cl^{\a, *}(M, E) $$
in the following way.
\begin{lem}\label{lem:cutoff}
Let $U$ be an open subset of $\R^n$ and  let $\sigma\in
\CS^{\a,k}(U,V)$ be a log-polyhomogeneous symbol of order $\a$ and
log-degree $k$. Then for any $x\in U$ the integral
$\int_{B_x^*(0,
R)} \sigma(x,\xi) \,\dbar\xi$ has an asymptotic expansion as
$R\to\infty$ $$\int_{B_x^*(0, R)} \sigma(x,\xi) \,\dbar\xi \
\sim_{ R\to \infty} $$
\begin{equation}\label{asympR}
C_x(\sigma) \ +\ \sum_{j=0,\alpha-j+n\neq 0}^\infty \sum_{l=0}^k
P_l(\sigma_{\alpha-j, l})(\log R) \, R^{\alpha-j+n} + \sum_{l=0}^k
\frac{1}{l+1} \int_{S_x^*U} \sigma_{-n,l } (x,\xi) \, \dbar_S\xi
\, \log^{l+1}  R
\end{equation}
where $P_l(\sigma_{\alpha-j, l})(X)$ is a polynomial of degree $l$
with coefficients depending on $\sigma_{\alpha-j, l}$. Here
$B_x^*(0, R)$ stands for  the ball of radius $R$ in the cotangent
space $T_x^*M$  and $S_x^*U$ the unit sphere in the cotangent
space $T_x^*U$.
\end{lem}

\vskip 1mm

Discarding the divergences, we can therefore  extract a finite
part from the asymptotic expansion of $\int_{B(0, R)}\sigma(x,\xi)
\dbar\xi$:
\begin{defn}
The finite-part integral \footnote{This concept is closely related
to {\it partie finie} of Hadamard \cite{Ha}, hence the terminology
used here. However in the physics literature this is also known as
``cut-off regularization''.} of $\sigma\in \CS^{\alpha, k}(U,V)$
is defined to be the constant term in the asymptotic expansion
\eqref{asympR}
\begin{equation}\label{cutoffint0}
\cutoffint_{T_x^*U} \sigma(x,\xi) \, \dbar\xi := {\rm
LIM}_{R\to\infty } \int_{B_x^*(0,R)} \sigma(x,\xi) \, \dbar\xi :=
C_x(\sigma).
\end{equation}
\end{defn}

\noi The proof of the following formula \cite{Gr1},\cite{Pa} and
of \lemref{lem:cutoff} is included in Appendix B.

\begin{lem}\label{cutoffint}
For  $\sigma\in \CS^{\alpha, k}(U,V)$
\begin{eqnarray}\label{eq:cutoffint}
\cutoffint_{T_x^*U} \sigma(x,\xi)\, \dbar\xi &= & \sum_{j=0}^{N}
\int_{B_x^*(0, 1)} \sigma_{\alpha-j}(x, \xi) \, \dbar\xi \ + \
\int_{T_x^*U}\sigma_{(N)}(x, \xi) \, \dbar\xi \nonumber\\ & & +
\sum_{j=0, \alpha-j+n\neq 0}^{N}\sum_{l=0}^k
\frac{(-1)^{l+1}l!}{(\alpha-j+n)^{l+1}}\int_{S_x^*U}
\sigma_{\alpha-j,l}(x, \xi)\, \dbar_S\xi.\hskip 15mm
\end{eqnarray}
It is independent of $N> \re(\alpha)+n-1$.
\end{lem}

\noi The residue terms on the right-side of \eqref{eq:cutoffint}
measure anomalous behaviour in the finite-part integral.
Specifically, \eqref{eq:cutoffint} implies that for a rescaling
$R\to \mu R$
\begin{eqnarray}
{\rm LIM}_{R\to\infty}\int_{B_x^*(0,\mu \, R)}\sigma(x,\xi)
\,\dbar\xi
&=& {\rm LIM}_{R\to\infty} \int_{B_x^*(0, R)} \sigma(x,\xi)\, \dbar\xi\nonumber\\
& & + \  \sum_{l=0}^k \frac{\log^{l+1}\mu}{l+1} \,\int_{S_x^*U}
\sigma_{\alpha-j,l}(x, \xi)\, \dbar_S\xi \label{rescale}
\end{eqnarray}
(cf. Appendix B) and hence that the finite-part integral is
independent of a rescaling if $\int_{S_x^*U} \sigma_{-n,l }
(x,\xi) \ \dbar_S\xi$ vanishes for each integer $0\leq l\leq k$.
More generally, just as ordinary integrals obey the transformation
rule $\vert {\rm det} C\vert\cdot \int_{\R^n} f(C\xi) \,\dbar\xi=
\int_{\R^n} f(\xi) \,\dbar\xi$, \, one hopes for a similar
transformation rule for the regularized integral $
\cutoffint_{\R^n} \sigma(\xi) \,\dbar\xi$ when $\sigma$ is a
log-polyhomogeneous symbol in order to obtain a globally defined
density on $M$. That, however, is generally not the case in the
presence of a residue  as  the following theorem shows.
\begin{prop}\cite{Le} \label{prop:noninv}
The finite part integral of $\sigma\in \CS^{*,k}(U)$ is generally
not invariant under a transformation $C\in Gl_n(T_x^*U)$. One has,
\begin{eqnarray}
\vert{\rm det} C\vert\cdot  \cutoffint_{T_x^*U} \sigma(x,C\xi) \,
\dbar\xi &=& \cutoffint_{T_x^*U} \sigma(x,\xi)\, \dbar\xi \\
& & + \sum_{l=0}^k \frac{(-1)^{l+1}}{l+1} \int_{S_x^*U}\sigma_{-n,
l}(x,\xi) \log^{l+1} \vert C^{-1}\xi\vert \, \dbar\xi. \nonumber
\end{eqnarray}
\end{prop}
\begin{proof}
We refer the reader to the proof of Proposition 5.2 in \cite{Le}.
\end{proof}
\noi As a consequence, whenever $\int_{S_x^*U}\sigma_{-n,
l}(x,\xi) \log^{l+1} \vert C^{-1}\xi\vert \, \dbar\xi$ vanishes
for each integer $0\leq l\leq k$ and $x\in U$ one then recovers
the usual transformation property
\begin{equation}\label{e:TRdensity}
\vert {\rm det} C\vert\cdot  \cutoffint_{T_x^*U} \sigma(x,C\xi)
\,\dbar\xi=  \cutoffint_{T_x^*U} \sigma(x,\xi) \,\dbar\xi.
\end{equation}

With respect to a trivialization $E_{|U} \cong U \times V$, a
localization of $A\in\Cl^{\alpha, k}(M, E)$ in $\Cl^{\alpha, k}(U,
V)$ can be written
$$Af(x) = \int_{\R^n}\int_U e^{i(x-y).\xi}\,\textsf{a}(x,y,\xi)\,f(y) \, dy\,\dbar\xi$$
with amplitude $\textsf{a}\in \CS^{\a}(U\times U,V)$. Then with
$$\s_A(x,\xi) := \textsf{a}(x,x,\xi)\in \CS^{\a}(U,V)$$
we define
$$\TRx(A) \, dx :=
\cutoffint_{T_x^*M} {\rm tr}_x(\sigma_A(x,\xi)) \, \dbar\xi \,
dx.$$ If \eqref{e:TRdensity} holds for $\s = \sigma_A$ in each
localization it follows that $\TRx(A) \, dx$ is independent of the
choice of local coordinates. This is known in the following cases.
\begin{prop}\label{p:TRdefined}
Let $A\in \Cl^{\alpha, k}(M, E)$. In each of the following cases
$\TRx(A)\,dx$ defines an element of $\Ci(M,|\Omega|)$, that is, a
global density  on $M$:
\begin{enumerate}
    \item $\a \notin [-n,\o) \, \cap \, \Z$,
    \item $A$ (of integer order) is even-even and $M$ is
    odd-dimensional,
    \item $A$ (of integer order) is even-odd and $M$ is
    even-dimensional.
\end{enumerate}
\end{prop}
Cases (1) and (2) were shown in \cite{KV}, where the canonical
trace was first introduced, in terms of homogeneous distributions.
Case (1) was reformulated in \cite{Le} in terms of finite part
integrals and extended to log-polyhomogeneous symbols $k\geq 0$.
Case (3) was introduced in \cite{Gr1} where it was shown that (2)
and (3) may be included in the finite part integral formulation.
We refer there for details. Notice though that it is easily seen
that the integrals $\int_{S_x^*U}\sigma_{-n, l}(x,\xi) \log^{l+1}
\vert C^{-1}\xi\vert \, \dbar\xi$ vanish in each case; for (1),
there is no homogeneous component of the symbol of degree $-n$ and
so the integrals vanish trivially, while setting
$g(\xi):=\sigma_{-n, l}(x,\xi) \log^{l+1} \vert C^{-1}\xi\vert$,
for cases (2) and (3) one has $g(-\,\xi)=-g(\xi)$ and so the
vanishing is immediate by symmetry.

\begin{defn}\label{def:TR}
For a $\pdo$ $A\in \Cl^{\alpha, *}(M, E)$ satisfying one of the
criteria {\rm (1),(2),(3)} in \propref{p:TRdefined} the canonical
trace is defined by
$${\rm TR}(A):= \int_M dx\, \TRx(A).$$
\end{defn}

The case of $\pdo$s of non-integer order is all that is needed for
the general formulae we prove here, cases (2) and (3) of
\propref{p:TRdefined} will be relevant only for applications and
refinements. Case (2) in particular includes differential
operators on odd-dimensional manifolds, though this holds by
default in so far as $\TR$ vanishes on differential operators in
any dimension (noted also in \cite{Gr1}):
\begin{prop}\label{prop:cutoffdiff}   Let $A\in \Cl(M, E)$ be a {\rm
differential} operator with local symbol $\sigma_A$, then for any
$x\in M$
$$\TRx(A) :=  \cutoffint_{T_x^* M} \tr_x(\sigma_A(x,\xi)) \ \dbar\xi = 0.$$
\end{prop}
\begin{proof}  Since
$A$ is a differential operator, $\sigma_A(x,\xi)=
\sum_{|k|=0}^{\ord A} \sigma_{k}(x,\xi) $ with $k = (k_1,\ldots,
k_n)$ a multi-index with $k_i\in\Nf$ and
$\sigma_{k}(x,\xi)=a_k(x)\xi^{k}$ positively homogeneous (with the
previous notations we have $\sigma_{(N)}=0$ provided $N\geq \ord
A$). Its finite-part integral on the cotangent space at $x\in M$
therefore reads
\begin{eqnarray*}
\cutoffint_{T_x^*M}\sigma_A(x,\xi) \dbar\xi &=& \sum_{|k|=0}^{\ord
A} a_k(x)\, {\rm LIM}_{R\to
\infty} \int_{B_x^*(0, R)}   \xi^k  \dbar\xi\\
&=& \sum_{|k|=0}^{\ord A} a_k(x) \, {\rm LIM}_{R\to \infty} \left(
\int_0^R r^{|k|+n-1} dr\right) \int_{S_x^*M}\xi^k \dbar\xi
\end{eqnarray*}
which vanishes since ${\rm LIM}_{R\to \infty}  \frac{ R^{|k|+n
}}{|k|+n}
 = 0$.
\end{proof}
On commutators the canonical trace has the following more
substantial vanishing properties \cite{KV}, \cite{MN}, \cite{Le},
\cite{Gr1}, providing some justification for its name.
\begin{prop}\label{p:TRtrace}
Let $ A\in \Cl^{a,k}(M, E)$, $B\in \Cl^{b,l}(M, E)$. In each of
the cases
\begin{enumerate}
    \item $\a + \b \notin [-n,\o) \, \cap \, \Z$,
    \item $A$ and $B$ are both even-even or are both even-odd and $M$ is
    odd-dimensional,
    \item $A$ is even-even, $B$ is even-odd, and $M$ is
    even-dimensional,
\end{enumerate}
the canonical trace is then defined on the commutator $[A,B]$ and
is equal to zero,
$$\TR([A,B])=0.$$
\end{prop}

\vskip 3mm

 The canonical trace extends the usual operator trace defined on
the subalgebra $\Cl^{\ord<-n}(M,E)$ of $\pdo$s of real order
$\re(\a) < -\,n$, in so far as for $\pdo$s with (real) order less
than $-n$ finite part integrals coincide with ordinary integrals.
More precisely, if $K_A(x, y)$ denotes the Schwartz kernel of
$A\in\Cl^{\ord<-n}(M,E)$ in a given localization, then
$\sigma_A(x,\xi)$ is integrable in $\xi$ and  $K_A(x, x)\, dx=
\left( \int_{T_x^*M} \sigma_A(x,\xi) \,\dbar\xi \right)\, dx$
determines a global density on $M$, and one has
$${\rm tr}(A)= \int_M dx\,{\rm tr}_x\left( K_A(x,
x)\right)= \int_M dx\,\cutoffint_{T_x^*M}\tr_x( \sigma_A(x,\xi))\,
\dbar\xi = \TR(A) .$$

\subsection{Holomorphic families of symbols}

We consider next families of symbols depending holomorphically on
a complex parameter $z$. The definition is somewhat more delicate
than that used in \cite{KV} (or \cite{Le}) since  growth
conditions are imposed on each $z$-derivative of the symbol. This
is in order to maintain control of the full Laurent expansion.

First, the meaning here of holomorphic dependence on a parameter
is  as follows. Let $W\subset\C$ be a complex domain, let $Y$ be
an open subset of $\R^m$, and let $V$ be a vector space. A
function $p(z,\eta)\in\Ci(W\times Y,\End(V))$ is holomorphic at
$z_0\in W$ if for fixed $\eta$ with
$$p^{(k)}(z_0,\eta) = \dd_z^k (p(z,\eta))|_{z=z_0}$$
there is a  Taylor expansion in a neighbourhood $N_{z_0}$ of $z_0$
\begin{equation}\label{e:Texpansion}
p(z,\eta) = \sum_{k=0}^{\o}p^{(k)}(z_0,\eta)\,\frac{(z-z_0)^k}{k!}
\end{equation}
which is convergent, uniformly on compact subsets of $N_{z_0}$,
with respect to the (metrizable) topology on $\Ci(W\times
Y,\End(V))$ associated with the family of semi-norms
\begin{equation}\label{e:seminorms}
\|q\|_{m, K_1, K_2}  = \sup_{\stackrel{(z,\eta)\in K_1 \times
K_2}{r + |\mu|\leq m}} |\dd_z^r \dd_{\eta}^{\mu} q(z,\eta)|
\end{equation}
defined for $m\in\N$ and compact subsets $K_1\subset W$,
$K_2\subset\R^m$.

\begin{defn}\label{defn:holfamilies} Let $m$ be a non-negative integer, let $U$ be an open
subset of $\R^n$, and let $W$ be a domain in $\C$. A holomorphic
family of log-polyhomogeneous symbols parametrized by $W$ of order
$\a\in\Ci(W,\C)$ and of log-degree $m$  means a function
$$\s(z)(x,\xi) := \s(z,x,\xi)\in\Ci(W \times U \times \R^n,
\End V)$$ for which:
\begin{enumerate}
\item $\s(z)(x,\xi)$ is holomorphic at $z\in W$ as an element of
$\Ci(W \times U \times \R^n, \End V)$ and
\begin{equation}\label{e:logclassical}
 \s(z)(x,\xi) \sim \sum_{j\geq 0}
 \s(z)_{\a(z)-j}(x,\xi) \ \in\CS^{\a(z),m}(U,V),
\end{equation}
where the function $\a:W\to \C$ is holomorphic;

\item for any integer $N\geq 1$ the remainder
$$\sigma_{(N)}(z)(x, \xi):= \sigma(z)(x,\xi)- \sum_{j=0}^{N-1}
\sigma_{\alpha(z)-j}(z)(x, \xi)$$ is holomorphic in $z\in W$ as an
element of $\Ci(W \times U \times \R^n, \End V)$ with $k^{{\rm
th}}$ $z$-derivative
\begin{equation}\label{e:kthderivlogclassical}
\sigma^{(k)}_{(N)}(z)(x, \xi) := \dd_z^k(\sigma_{(N)}(z)(x, \xi))
\in {\rm S}^{\a(z)-N + \e}(U,V)
\end{equation}
for any $\e>0$.
\end{enumerate}

\vskip 2mm A  family $z\mapsto A(z)$ of log-classical $\pdo$s on
$\Ci(M,E)$ parametrized by a domain $W\subset \Cf$ is holomorphic
if in each local trivialisation of $E$ one has $$A(z) = {\rm
Op}(\sigma_{A(z)}) + R(z)$$ with $\sigma_{A(z)}$  a holomorphic
family of log-polyhomogeneous symbols and $R(z)$ a smoothing
operator with Schwartz kernel $R(z,x,y)\in \Ci(W\times X\times
X,\End(V))$ holomorphic in $z$.
\end{defn}
There are, of course, other ways to express these conditions; for
example, in terms of the truncated kernel $K^{(N)}(z)(x, y):=
\int_{T_x^*U} e^{i \xi\cdot (x-y)} \sigma_{(N)}(z)(x, \xi)\,
\dbar\xi$ with large $N$, and its derivatives $\dd^k_z
K^{(N)}(z)(x, y)$, used in the case $k=0$ in \cite{KV}  to compute
the pole of $\Tr(A(z))$ at $z_0\in P$. When dealing with the full
Laurent expansion the essential requirement is that a balance be
preserved between the Taylor expansion in $z$, in terms of the
growth rates of the $z$-derivatives of the symbol, and the
asymptotic symbol expansion in $\xi$.

\vskip 1mm

\begin{prop}\label{p:asympsymbolderiv}
If $\s(z)(x,\xi) \in\CS^{\a(z),m}(U,V)$ is a holomorphic
family of log-classical symbols, then so is each derivative
\begin{equation}\label{derivsigma}
\s^{(k)}(z)(x,\xi) := \dd_z^k\left(\sigma(z)(x, \xi)\right) \
\in\CS^{\a(z),m+k}(U,V).
\end{equation}
Precisely, $\s^{(k)}(z)(x,\xi) $ has asymptotic expansion
\begin{equation}\label{asympderivsigma}
\s^{(k)}(z)(x,\xi) \sim \sum_{j\geq 0}
\s^{(k)}(z)_{\a(z)-j}(x,\xi)
\end{equation}
where as elements of $\bigcup_{l=0}^{m+k} \CS^{\a(z)-j,\,l}(U,V)$
\begin{equation}\label{componentsderivsigma}
\s^{(k)}(z)_{\a(z)-j}(x,\xi) = \dd_z^k\left(\sigma(z)_{\a(z)-j}(x,
\xi)\right).
\end{equation}
That is,
\begin{equation}\label{componentsderivsigma2}
\left(\dd_z^k\sigma(z)\right)_{\a(z)-j}(x,\xi) =
\dd_z^k\left(\sigma(z)_{\a(z)-j}(x,
   \xi)\right).
\end{equation}
\end{prop}
\begin{proof} We have to show that
\begin{equation}\label{derivsigma2}
\dd_z^k\left(\sigma(z)(x, \xi)\right) \sim \sum_{j\geq 0}
\dd_z^k\left(\sigma(z)_{\a(z)-j}(x, \xi)\right)
\end{equation}
where the summands are log-polyhomogeneous of the asserted order.
First, the estimate
$$\dd_z^k\left(\sigma(z)(x, \xi)\right) - \sum_{j=0}^{N-1}
\dd_z^k\left(\sigma(z)_{\a(z)-j}(x, \xi)\right) \in {\rm
S}^{\a(z)- N + \e}(U,V)$$ any $\e>0$, needed for
\eqref{derivsigma2} to hold is equation
\eqref{e:kthderivlogclassical} of the definition. It remains to
examine the form of the summands in $\sum_{j=0}^{N-1}
\dd_z^k\left(\sigma(z)_{\alpha(z)-j}(x, \xi)\right).$ Taking
differences of remainders $\sigma_{(N)}(z)(x, \xi)$ implies that
each term $\sigma(z)_{\a(z)-j}(x, \xi)$ is holomorphic. In order
to compute $\dd_z\left(\sigma(z)_{\alpha(z)-j}(x, \xi)\right)$
one must compute the derivative of each of its homogeneous
components; for $|\xi|\geq 1$ and any $l\in \{0, \cdots, m\}$
\begin{eqnarray*}
\dd_z\left( \sigma_{ \alpha(z) -j,\, l} (z)(x, \xi)\right)&=&
\dd_z\left( \vert \xi\vert^{\alpha(z)-j}\ \sigma_{ \alpha(z) -j,\,
l} (z)(x, \frac{\xi}{\vert\xi\vert})\right)\\ &=&
\left(\alpha^\prime(z)\, \vert \xi\vert^{\alpha(z)-j}\ \sigma_{
\alpha(z) -j,\,l} (z)(x, \frac{\xi}{\vert\xi\vert})\right) \log
\vert \xi\vert \\ & & + \ \ \vert \xi\vert^{\alpha(z)-j}\
\dd_z\left( \sigma_{ \alpha(z) -j,\,l} (z)(x,
\frac{\xi}{\vert\xi\vert})\right).
\end{eqnarray*}
Since $ \sigma_{ \alpha(z) -j,\,l}(z)(x, \xi\vert\xi\vert\ii)$ is
a symbol of constant order zero, so is its $z$-derivative. Hence
\begin{equation}\label{e:derivsigma}
\dd_z\left(\sigma(z)_{\alpha(z)-j,\,l}(x, \xi)\right) =
\alpha^\prime(z)\, \sigma(z)_{ \alpha(z) -j,\,l}(x, \xi)\,\log
[\xi] \ + \ p_{ \alpha(z) -j,\,l}(z)(x, \xi)
\end{equation}
where $\sigma_{ \alpha(z) -j,\,l}(z),\, p_{ \alpha(z)
-j,\,l}(z)\in \CS^{\a(z)-j}(U)$ are homogeneous in $\xi$ of order
$\alpha(z) -j$. Hence
$\dd_z\left(\sigma(z)_{\alpha(z)-j}\right)\in
\CS^{\a(z)-j,\,m+1}(U)$. Iterating \eqref{e:derivsigma},
$\dd_z^k\left(\sigma_{\alpha(z)-j}(z)(x, \xi)\right)$ is thus seen
to be a polynomial in $\log[\xi]$ of the form
$$(\alpha^\prime(z))^k\,
\sigma_{ \alpha(z) -j,\,k}(z)(x, \xi)\,\log^{k+m} [\xi]\ + \
\ldots \ + \  \ \vert \xi\vert^{\alpha(z)-j}\ \dd^k_z( \sigma_{
\alpha(z) -j,\,l} (z)(x, \frac{\xi}{\vert\xi\vert})) \,\log^{0}
[\xi]
$$ with each coefficient homogeneous
of order $\alpha(z) -j$. This completes the proof.
\end{proof}
\vskip 1mm

Thus, taking derivatives adds more logarithmic terms to each term
$\s(z)_{\a(z)-j}(x,\xi)$, increasing the log-degree, but the order
is unchanged. Specifically, $\s^{(k)}(z)_{\a(z)-j}$ takes the form
\begin{equation}\label{componentsderivsigma3}
\s^{(k)}(z)_{\a(z)-j}(x,\xi) =  \sum_{l=0}^{m+k}
\s^{(k)}(z)_{\a(z)-j,\,l}(x,\xi)\log^l[\xi],
\end{equation}
where the terms $\s^{(k)}(z)_{\a(z)-j,\,l}(x,\xi)$ are positively
homogeneous in $\xi$ of degree $\a(z)-j$ for $|\xi|\geq 1$ and can
be computed explicitly from the lower order derivatives of
$\sigma(z)_{\a(z)-j,\,m}(x, \xi)$. The following more precise
inductive formulae will be needed in what follows.
\begin{lem}\label{inductivederivs}
Let $\s(z)(x,\xi)\in\CS(U,V)$ be a holomorphic family of classical
symbols. Then for $|\xi|\geq 1$
\begin{eqnarray*}
\sigma^{(k+1)}_{\alpha(z)-j,\,k+1}(z)(x, \xi) &= & \a^\prime(z)\,
\sigma^{(k)}_{\alpha(z)-j,\,k}(z)(x, \xi),\\[2mm]
\sigma^{(k+1)}_{\alpha(z)-j,\,l}(z)(x, \xi) &= & \a^\prime(z)\,
\sigma^{(k)}_{\alpha(z)-j,\,l-1}(z)(x, \xi) \\[2mm] & & + \
|\xi|^{\alpha(z) -j}\,\dd_z(\sigma^{(k)}_{ \alpha(z) -j,\,l}
(z)(x, \xi/|\xi|)), \ \ \ \
1\leq l\leq k,\\[2mm]
\sigma^{(k+1)}_{\alpha(z)-j,\,0}(z)(x, \xi) & = & |\xi|^{\alpha(z)
-j}\,\dd_z(\sigma^{(k)}_{ \alpha(z) -j,0} (z)(x, \,\xi/|\xi|)).
\end{eqnarray*}
\end{lem}
\begin{proof}
From the above
\begin{equation*}
\s^{(k)}(z)_{\a(z)-j}(x,\xi) =   \dd^k_z(\s(z)_{\a(z)-j}(x,\xi)) =
\sum_{l=0}^{k} \s^{(k)}(z)_{\a(z)-j,\,l}(x,\xi)\log^l[\xi],
\end{equation*}
so that
\begin{equation}\label{sigmakplusone}
\s^{(k+1)}(z)_{\a(z)-j}(x,\xi) = \sum_{l=0}^{k}
\dd_z\left(\s^{(k)}(z)_{\a(z)-j,\,l}(x,\xi)\right)\log^l[\xi].
\end{equation}
Hence for $|\xi|\geq 1$
$$\sum_{l=0}^{k+1} \s^{(k+1)}(z)_{\a(z)-j,\,l}(x,\xi)\log^l|\xi|
=$$
\begin{equation*}
\sum_{r=0}^{k}\alpha^\prime(z)\, \sigma^{(k)}(z)_{ \alpha(z)
-j,\,r}(x, \xi)\,\log^{r+1} |\xi| \ + \ |\xi|^{\alpha(z)-j}\
\dd_z\left( \sigma^{(k)}(z)_{ \alpha(z) -j,\,r} (x,
\frac{\xi}{|\xi|})\right)\, \log^r |\xi|
\end{equation*}
where for the right-side we apply \eqref{e:derivsigma} to each of
coefficient on the right-side of \eqref{sigmakplusone}.  Equating
coefficients completes the proof.
\end{proof}

A corresponding result on the level of operators follows in a
straightforward manner:
\begin{prop}\label{prop:diffsymb}
Let $z\mapsto A(z)\in \Cl^{\a(z), m}(M, E)$ be a holomorphic
family of log-polyhomogeneous  $\pdo$s. Then for any non-negative
integer $k$, $A^{(k)}(z_0)$ lies in $\Cl^{\a(z_0), m+k}(M, E)$.
\end{prop}

\vskip 2mm
\begin{ex}\label{example}{\rm
For real numbers $\a, q$ with $q>0$ the function
$\sigma(z)(x,\xi)=\psi(\xi)\,\vert \xi \vert^{\a- q z}$, where
$\psi$ is a smooth cut-off function which vanishes near the origin
and is equal to $1$ outside the unit ball, provides a holomorphic
family of classical symbols; at any point $z=z_0\in \C$ we have
$\sigma^{(k)}(z_0)(x, \xi)=(-q)^k\psi(\xi)\,\log^k \vert
\xi\vert\, \vert \xi\vert^{\a-q z_0}$ which lies in $\CS^{\a-q
z_0, k}(U)$. More generally, if $Q\in \Cl^{q}(M,E)$ is a classical
elliptic $\pdo$ of order $q> 0$ with principal angle $\th$, then
one has for each $z \in \Cf$ the complex power $Q_{\th}^{-z}\in
\Cl^{-qz}(M,E)$ \cite{Se1} represented in a local coordinate chart
$U$ by a classical symbol ${\bf q}(z)(x,\xi)\in\CS^{-qz}(U,V)$.
Let $A\in \Cl^{\a}(M,E)$ be a coefficient classical $\pdo$
represented in $U$ by ${\bf a}(x,\xi)\in\CS^{\a}(U)$. Then
$\s_{AQ^{-z}_{\th}}(x,\xi)\in\CS^{\a-qz}(U,V)$ is a holomorphic
family of symbols parametrized by $W=\C$ whose convergent Taylor
expansion in $\Ci(\C\times U,V)$ around each $z_0\in\C$ is from
\cite{O1} Lemma(2.1) given by
$$(\s_{AQ^{-z}_{\th}})_{\a-qz-j}(x,\xi) =
\sum_{k=0}^{\o}\sum_{l=0}^k(-1)^k\left({\bf a}\circ \log^k({\bf
q})\circ{\bf q}(z_0)\right)_{\a-qz_0 -j,\,l}(x,\xi)
\,\log^l|\xi|\, \frac{(z-z_0)^k}{k!},$$ where $\circ$ denotes the
usual mod(${\rm S}^{-\o})$ symbol product, ${\bf q} := {\bf
q}(-1)$ and $\log^k({\bf q})(x,\xi):= (\log({\bf  q})\circ
\ldots\circ \log({\bf q}))(x,\xi)\in\CS^{0,k}(U,V)$ with $k$
factors.}
\end{ex}
\vskip 2mm

\subsection{A Laurent expansion for finite part integrals of holomorphic symbols}

The following theorem computes the Laurent expansion for finite
part integrals of holomorphic families of classical symbols of
order $\alpha(z)$ in terms of local canonical and residue
densities. This extends Proposition 3.4 in \cite{KV}, and results
of \cite{Gu,Wo}, where the pole, the first coefficient in the
expansion, was identified as the residue trace. The proof uses the
property that each term of the Taylor series of a holomorphic
family of classical symbols has an asymptotic symbol expansion,
allowing the Laurent expansion of $\cutoffint \,\sigma(z)( x,\xi)
\, \dbar\xi$ to be computed through \lemref{cutoffint}. Notice
that although the Taylor expansion in the $\Ci$ topology gives no
control over the symbol as $|\xi| \to\o$, \eqref{e:logclassical},
\eqref{e:kthderivlogclassical} impose what is needed to ensure
integrability requirements.

\begin{defn}
A holomorphic function $\a:W\to\Cf$ defined on a domain $W\subset
\Cf$ is said to be  non-critical on
$$P:=\alpha^{-1}(\Z\cap [-n, +\infty[)\cap W$$ if
$\alpha^\prime(z_0)\neq 0$ at each $z_0\in P$.
\end{defn}

\vskip 3mm

\begin{thm}\label{thm:KV2}
\noi {\rm (1)} \ Let $U$ be an open subset of $\R^n$. Let
$z\mapsto\sigma(z)\in \CS^{\a(z)}(U,V)$ be a holomorphic family of
classical symbols parametrized by a domain  $W\subset \C$  such
that the order function $\a$ is non-critical on $P$. Then for each
$x\in U$ the map $z\mapsto \cutoffint_{T_x^*U} \sigma(z)( x,\xi)
\, \dbar\xi$ is a meromorphic function on $W$  with poles located
in $ P$. The poles are at most simple and for $z$ near $z_0\in P$
one has

\begin{eqnarray} \label{eq:KV2symb}
& & \cutoffint_{T_x^*U} \sigma(z)(x, \xi) \, \dbar\xi \  =  \ - \
\frac{1}{ \alpha^\prime( z_0 )}\,\int_{S_x^*U} \sigma(z_0)_{-n}
(x,\xi) \ \dbar_S\xi\ \frac{1}{
(z-z_0)}\nonumber\\[3mm]
& & \  + \ \left( \cutoffint_{T_x^*U}\, \sigma(z_0) (x,
\xi)\,\dbar\xi \ - \ \frac{1} {\alpha^\prime(z_0)}\int_{S_x^*U}
\sigma^\prime(z_0)_{-n,0} (x,\xi) \ \dbar_S\xi\right) \nonumber\\[3mm]
& & \hskip 43mm \, \ \ + \ \
\frac{\a^{\prime\prime}(z_0)}{2\,\a^\prime(z_0)^2}\ \int_{S_x^*U}
\sigma(z_0)_{-n} (x,\xi) \ \dbar_S\xi  \nonumber\\[1mm]
 &  & \  + \ \sum_{k=1}^K
\left(\cutoffint_{T_x^*U}\sigma^{(k)}(z_0)(x, \xi)\,
\dbar\xi\newline - \int_{S_x^*U}{\mathcal L}_k(\sigma(z_0),
\cdots,
\sigma^{(k+1)}(z_0))_{-n,0} (x,\xi) \ \dbar_S\xi\,\right)\frac{(z-z_0)^k}{k!}\nonumber\\[3mm]
\ \ & &   \ + \ \ o\left((z-z_0)^K\right),
\end{eqnarray}

\vskip 2mm

\noi where
\begin{equation}\label{e:L}
{\mathcal L}_k\left(\sigma(z_0), \cdots,
\sigma^{(k+1)}(z_0)\right) =
\sum_{j=0}^{k+1}\frac{p_{k+1-j}}{\small{\a^{'}(z_0)}^{k+2-j}}\
\sigma^{(j)}(z_0) \ \ \in \CS^{\a(z_0),k+1}(U,V),
\end{equation}
and $p_{k+1-j}$ is an explicitly computable polynomial of degree
$k+1-j$ in \ $\a^{'}(z_0),\ldots,\a^{(k+1)}(z_0)$. Furthermore,
the coefficient of $\frac{(z-z_0)^k}{k!}$ in \eqref{eq:KV2symb} is
equal to ${\rm fp}_{z=z_0} \cutoffint_{T_x^*U}\sigma^{(k)}(z).$ If
$\a$ is a linear function $\alpha(z) = qz + b$ with $q\neq 0$ then
\eqref{eq:KV2symb} reduces to
\begin{eqnarray}
& & \cutoffint_{T_x^*U} \sigma(z)(x, \xi) \, \dbar\xi \  =   \ \ -
\ \ \frac{1}{q}\,\int_{S_x^*U} \sigma(z_0)_{-n} (x,\xi) \
\dbar_S\xi\ \frac{1}{
(z-z_0)}\nonumber\\[3mm]
& & \ \ \ +\ \ \left( \cutoffint_{T_x^*U}\, \sigma(z_0) (x,
\xi)\,\dbar\xi \ - \ \frac{1} {q}\int_{S_x^*U}
\sigma^\prime(z_0)_{-n,0} (x,\xi) \ \dbar_S\xi\right) \nonumber\\[3mm]
 &  & \ \ \ +\ \ \sum_{k=1}^K
\left(\cutoffint_{T_x^*U}\sigma^{(k)}(z_0)(x, \xi)\,
\dbar\xi\newline \ -  \ \frac{1}{q(k+1)}\,\int_{S_x^*U}
\sigma^{(k+1)}(z_0)_{-n,0} (x,\xi) \ \dbar_S\xi\right)\, \frac{(z-z_0)^k}{k!}\nonumber\\
[3mm] & &   \ \ \ + \ \ o\left((z-z_0)^K\right).
\label{eq:KV2symbq}
\end{eqnarray}

\vskip 2mm

\noi If $z_0\in W$ but $z_0\notin P$, then $\cutoffint_{T_x^*U}
\sigma(z)(x, \xi) \, \dbar\xi $ is holomorphic at $z=z_0$ and
\eqref{eq:KV2symb} then simplifies to the Taylor expansion
$$\cutoffint_{T_x^*U} \sigma(z)(x, \xi) \, \dbar\xi \  = \
\cutoffint_{T_x^*U}\, \sigma(z_0) (x, \xi)\,\dbar\xi  \  + \
\sum_{k=1}^K \cutoffint_{T_x^*U}\sigma^{(k)}(z_0)(x, \xi)\,
\dbar\xi \, \frac{(z-z_0)^k}{k!}$$
\begin{equation}\label{eq:KV2symbhol}
+  \ o\left((z-z_0)^K\right).
\end{equation}
\vskip 3mm

\noi {\rm (2)} \  For any holomorphic family $z\mapsto A(z)\in
\Cl^{\a(z)}(M, E)$ of classical \, $\pdo$s \, parametrized by a
domain $W\subset \C$, such that order function $\a$ is
non-critical on $P$, the map $z\mapsto {\rm TR}(A(z)):=  \int_M \,
dx\int_{T^*_xM}\, {\rm tr}_{x} \left(\sigma_{A(z)} (x, \xi)
\right) \, \dbar\xi$ is a meromorphic function on $W$  with  poles
located in  $ P$. The poles are at most simple and for $z$ near
$z_0\in P$

\begin{eqnarray} \label{eq:KV2op}
\TR(A(z) )&=& -  \ \frac{1}{\alpha^\prime(z_0)} \, {\rm
res}(A(z_0))\,
\frac{1}{(z-z_0)} \nonumber\\[3mm] & & + \ \int_M dx\, \left(\TRx(A(z_0))
-\frac{1}{\alpha^\prime(z_0)}\, {\rm
res}_{x,0}(A^\prime(z_0))\right) \ \ + \ \
\frac{\a^{\prime\prime}(z_0)}{2\,\a^\prime(z_0)^2}\
\res\,\left(A(z_0)
\right)\nonumber\\[3mm]
& & +  \ \sum_{k=1}^K \int_M\,dx\, \left(\TRx(A^{(k)}(z_0)) -
\resxo\left({\mathcal L}_k(\sigma_{A(z_0)}, \cdots,
\sigma_{A^{(k+1)}(z_0)})\right)\right) \frac{(z-z_0)^k}{k!} \nonumber\\[3mm]
& & + \ o\left((z-z_0)^K\right).
\end{eqnarray}
\noi  Furthermore, the coefficient of $\frac{(z-z_0)^k}{k!}$ in
\eqref{eq:KV2op} is equal to ${\rm fp}_{z=z_0}\TR(A^{(k)}(z)).$ If
$A(z)$ has order $\alpha(z) = qz + b$ with $q\neq 0$ then
\begin{eqnarray}
\TR(A(z))&=& -  \ \frac{1}{q} \, {\rm res}(A(z_0))\,
\frac{1}{(z-z_0)} \nonumber\\[3mm] & & + \ \int_M dx\, \left(\TRx(A(z_0))
-\frac{1}{q}\,{\rm res}_{x,0}(A^\prime(z_0))\right) \nonumber\\
[3mm] & & +  \ \sum_{k=1}^K \int_M\,dx\, \left(\TRx(A^{(k)}(z_0))
- \frac{\resxo\left(\sigma_A^{(k+1)}(z_0)\right)}{
q\,(k+1)}\right)\, \frac{(z-z_0)^k}{k!}\nonumber\\
[3mm] & & + \ o\left((z-z_0)^K\right). \label{linearorder}
\end{eqnarray}
\vskip 2mm \noi If $z_0\in W$ but $z_0\notin P$, then $\TR(A(z))$
is holomorphic at $z=z_0$ and \eqref{eq:KV2op} then simplifies to
the Taylor expansion
\begin{equation}\label{eq:KV2ophol}
\TR(A(z)) \  = \ \TR(A(z_0))   +  \sum_{k=1}^K \TR(A^{(k)}(z_0))
\,\frac{(z-z_0)^k}{k!} +   o\left((z-z_0)^K\right).
\end{equation}
\end{thm}

\vskip 5mm

\begin{rk}
Since $\a$ is non-critical on $P$, we have from
\propref{prop:diffsymb} and equation \eqref{e:derivsigma} that the
operators $A^{(k)}(z_0)\in \Cl^{\a(z_0),k}(M,E)$ in equation
\eqref{eq:KV2op} are not classical for $k\geq 1$.
\end{rk}

\begin{rk} At a point $z_0\in P$, $\alpha^\prime(z_0)\neq 0$; writing
$\alpha(z)= \alpha(z_0)+ \alpha^\prime(z_0)(z-z_0) +o(z-z_0)$ we
find that $\alpha$ is injective in a neighborhood of $z_0$. As a
consequence, $\Z$ being countable, so is the set of poles
$P=\alpha^{-1}(\Z\cap [-n, +\infty)\,)\cap W$ countable.
\end{rk}

\begin{rk}\label{rk} Setting $\alpha(z)= \frac{z}{1+\lambda\, z}$ with $\lambda\in \R^*$
for $z\in\C\backslash\{-\la\ii\}$ gives rise to an additional
finite part $\frac{\alpha^{\prime\prime}(0)}{2\alpha^\prime(0)^ 2}
\int_{S_x^*U} \sigma(0)_{-n} (x,\xi) \ \dbar_S\xi= \lambda\,
\int_{S_x^*U} \sigma(0)_{-n} (x,\xi) \ \dbar_S\xi$ just as  a
rescaling $R\to e^{\lambda}\,R$ in the finite part integrals gives
rise to the  extra term $\lambda\, \int_{S_x^*U} \sigma(0)_{-n}
(x,\xi) \ \dbar_S\xi$ (see  \eqref{rescale} with $k=0$ and $\mu=
e^{\lambda}$).
\end{rk}

\begin{proof}
Since the  orders $\alpha(z)$  define a holomorphic map at each
point of  $P$, for any $z_0\in P$  there is a ball $B(z_0, r)
\subset W\subset \C$ centered at $z_0\in W$ with radius $r>0$ such
that $\left( B(z_0, r)\backslash\{z_0\}\right)\, \bigcap\,  P=
\phi $. In particular, for all $ z\in B(z_0, r)\backslash\{z_0\}$,
the symbols $\sigma(z)$ have non-integer order.  As a consequence,
outside the set $P$, the finite part integral $\cutoffint_{T_x^*U}
\sigma(z)( x,\xi)\dbar\xi$ is defined without ambiguity and
$\cutoffint_{T_x^*U} \sigma(z)( x,\xi)\dbar\xi\,dx$ defines a
global density on $M$.

Since $z_0 \in P$,  there is some $j_0\in \N\cup \{0\}$ such that
$\alpha(z_0)+n-j_0=0$. On the other hand, for $z\in B(z_0,
r)\backslash\{z_0\}$ we have $\a(z) + n -j\neq 0$ and
$N>\re(\alpha(z))+n-1$ can be chosen uniformly to ensure that
$\sigma_{(N)}(z) \in {\rm S}^{\,<-n}(U,V)$. Hence for $z\in B(z_0,
r)\backslash\{z_0\}$ equation (\ref{eq:cutoffint}) yields (with
$k=0$)
$$\cutoffint_{T_x^*U}
\sigma(z)(x, \xi)\ \dbar\xi$$
\begin{eqnarray}
& = &\sum_{j=0}^N \int_{B_x^*(0,1)} \sigma(z)_{\alpha(z)-j}(x,\xi)
\ \dbar\xi \ + \ \int_{T_x^*U} \sigma_{(N)}(z)(x,\xi) \ \dbar\xi
\nonumber\\ \hskip 10mm &
 & - \ \sum_{j=0 }^N \frac{1}{\alpha(z)+n-j} \int_{S_x^*U}
\sigma(z)_{\alpha(z)-j}(x,\xi) \, \dbar_S\xi\nonumber\\
&  =& \sum_{j=0}^N\int_{B_x^*(0,1)}\sigma(z)_{\alpha(z)-j}(x,\xi)
\ \dbar\xi +\int_{T_x^*U} \sigma_{(N)}(z)(x,\xi) \ \dbar\xi \nonumber\\
&& - \ \sum_{j=0 , j\neq j_0 }^N\frac{1}{\alpha(z)+n-j}
\int_{S_x^*U}\sigma(z)_{\alpha(z)-j}(x, \xi)\, \dbar_S\xi\,\nonumber\\
&& - \ \frac{1}{\alpha(z)-\alpha(z_0)}
\int_{S_x^*U}\sigma(z)_{\a(z)-j_0}(x, \xi)\, \dbar_S\xi\,
\label{e:cutoffsigmaz}
\end{eqnarray}
where, in view of the growth conditions \eqref{e:logclassical} and
\eqref{e:kthderivlogclassical}, it is not hard to see that each of
the integrals on the right-side of \eqref{e:cutoffsigmaz} is
holomorphic in $z$. Since $\sigma_{ \alpha(z) -j}(z)(x, \xi)$ is a
holomorphic family of classical symbols, there is a Taylor
expansion \eqref{e:Texpansion}
\begin{equation}\label{e:sigmazexpansion0}
\sigma(z)_{ \alpha(z) -j}(x, \xi) = \sum_{k=0}^{\o}\s^{(k)}(z_0)_{
\alpha(z_0) -j}(x, \xi)\,\frac{(z-z_0)^k}{k!}
\end{equation}
with coefficients in $\CS^{\a(z_0)-j,\, k}(U)$
\begin{equation}\label{e:coeffts}
\s^{(k)}(z_0)_{ \alpha(z_0) -j}(x,\xi) : = \dd^k_z \left(\s(z)_{
\alpha(z) -j}\right)|_{z=z_0} = \left(\dd^k_z
\s(z)\right)_{\alpha(z) -j} |_{z=z_0} ,
\end{equation}
where the first equality is by definition while the second
equality is equation \eqref{componentsderivsigma2}, and likewise
there is a Taylor expansion of the remainder
$\sigma_{(N)}(z)(x,\xi)$ with coefficients
\begin{equation}\label{e:Rcoeffts}
\dd^k_z \left(\sigma_{(N)}(z)(x,\xi)\right) |_{z=z_0} =
\left(\dd^k_z \sigma(z)\right)_{(N)}(x,\xi) |_{z=z_0},
\end{equation}
where again the equality is consequent on equations
\eqref{asympderivsigma} and \eqref{componentsderivsigma2}. For any
non-negative integer $K$ we may therefore rewrite the first two
lines of \eqref{e:cutoffsigmaz} as a polynomial $\sum_{k=0}^K\,a_k
\, \frac{(z-z_0)^k}{k!} $ plus an error term of order
$o((z-z_0)^K)$ with
\begin{eqnarray}
a_k  &  = & \sum_{j=0}^N\int_{B_x^*(0,1)}\left(\dd^k_z
\s(z)\right)_{ \alpha(z) -j}(x,\xi)|_{z=z_0}\ \dbar\xi
+\int_{T_x^*U} \left(\dd^k_z \sigma(z)\right)_{(N)}(x,\xi) |_{z=z_0} \ \dbar\xi \nonumber\\
& & -\sum_{j=0 , j\neq j_0 }^N
\left.\dd_z^k\right|_{z=z_0}\left(\frac{1}{\alpha(z)+n-j}
\int_{S_x^*U}\sigma(z)_{\alpha(z)-j}(x, \xi)\
\dbar_S\xi\,\right).\label{ak}
\end{eqnarray}
Here, since $j\neq j_0$, we use the fact that each factor in the
terms of the final summation of \eqref{ak} are holomorphic in a
neighbourhood of $z_0$ (including at $z=z_0$). On the other hand,
from \eqref{asympderivsigma}
$$
\s^{(k)}(z_0)(x,\xi) := \dd^k_z \s(z)(x,\xi)|_{z=z_0} \sim
\sum_{j\geq 0} \left(\dd^k_z
\s(z)\right)_{\a(z)-j}(x,\xi)|_{z=z_0}
$$
while we know from \eqref{derivsigma} that $\s^{(k)}(z) \in
\CS^{\a(z), k}(U)$. Hence \eqref{eq:cutoffint} may be applied to
see that
$$\cutoffint_{T_x^*U}
\s^{(k)}(z_0)(x,\xi) \ \dbar\xi$$
\begin{eqnarray}
=& &\sum_{j=0}^N \int_{B_x^*(0,1)} \left(\dd^k_z
\s(z)\right)_{\alpha(z) -j}(x,\xi)  |_{z=z_0}\ \dbar\xi \ + \
\int_{T_x^*U}
\left(\dd^k_z \sigma(z)\right)_{(N)}(x,\xi) |_{z=z_0}\  \dbar\xi \nonumber\\
\hskip 10mm & + & \  \sum_{j=0,j\neq j_0}^{N}\sum_{l=0}^k
\frac{(-1)^{l+1}l!}{(\alpha(z_0)-j+n)^{l+1}}\int_{S_x^*U}
\left(\dd_z^k\s(z)\right)_{\alpha(z)-j,\,l}(x,\xi)|_{z=z_0}\
\dbar_S\xi \label{e:cutoffsigmaderivative}
\end{eqnarray}

\noi From the following lemma we conclude  that the expressions in
\eqref{ak} and \eqref{e:cutoffsigmaderivative} are equal.
\begin{lem}\label{lem:transition}
For $j\neq j_0$ one has in a neighbourhood of $z_0$

$$\dd_z^k\left(\frac{-1}{\alpha(z)+n-j}
\int_{S_x^*U}\sigma(z)_{\alpha(z)-j}(x, \xi)\,
\dbar_S\xi\,\right)\hskip 20mm
$$
\begin{equation}\label{e:transition}
\hskip 40mm = \ \ \sum_{l=0}^k
\frac{(-1)^{l+1}l!}{(\alpha(z)-j+n)^{l+1}}\int_{S_x^*U}
\left(\dd_z^k\s(z)\right)_{\alpha(z)-j,\,l}(x,\xi)\, \dbar_S\xi.
\end{equation}

\end{lem}
\begin{proof}
We choose $z$ in a neighbourhood of $z_0$ such that each of the
factors on both sides of \eqref{e:transition} are holomorphic. The
equality holds trivially for $k=0$. For clarity we check the case
$k=1$ before proceeding to the general inductive step. For $k=1$
the left-side of \eqref{e:transition} is equal to
$$ \frac{\a^{\prime}(z)}{(\alpha(z)-j+n)^2}\int_{S_x^*U} \sigma(z)_{\alpha(z)-j}(x,
\xi)\, \dbar_S\xi \hskip 40mm$$
\begin{equation}\label{e:transition1}
 \hskip 40mm - \ \ \frac{1}{\alpha(z)-j+n}\int_{S_x^*U}
\dd_z\left(\s(z)_{\alpha(z)-j}\right)(x,\xi)\, \dbar_S\xi.
\end{equation}
From \eqref{asympderivsigma} and \eqref{e:derivsigma}, for
$|\xi|\geq 1$
\begin{equation*}
\left(\dd_z\sigma(z)\right)_{\alpha(z)-j}(x, \xi) =
\alpha^\prime(z)\, \sigma(z)_{ \alpha(z) -j}(x, \xi)\,\log |\xi| \
+ \ p_{ \alpha(z) -j}(z)(x, \xi)
\end{equation*}
and hence $\left(\dd_z\s(z)\right)_{\alpha(z)-j,\,1}(x,\xi) =
\alpha^\prime(z)\, \sigma(z)_{ \alpha(z) -j}(x, \xi)$ for
$|\xi|\geq 1$. The expression in \eqref{e:transition1} is
therefore equal to
$$ \frac{1}{(\alpha(z)-j+n)^2}\int_{S_x^*U} \left(\dd_z\s(z)\right)_{\alpha(z)-j,\,1}(x,\xi)
\, \dbar_S\xi \hskip 30mm$$ $$\hskip 12mm - \ \
\frac{1}{\alpha(z)-j+n}\int_{S_x^*U}
\dd_z\left(\s(z)_{\alpha(z)-j}\right)(x,\xi)\, \dbar_S\xi$$ which
is the right-side of \eqref{e:transition} for $k=1$.

Assume now that \eqref{e:transition} holds for some arbitrary
fixed $k\geq 0$. Then the left-side of \eqref{e:transition} for
$k+1$ is equal to
$$\dd_z\left(\dd_z^{k}\left({-1\over \alpha(z)+n-j}
\int_{S_x^*U}\sigma(z)_{\alpha(z)-j}(x, \xi)\,
\dbar_S\xi\,\right)\right)$$
\begin{eqnarray}
&  = & \dd_z\left(\sum_{l=0}^k
\frac{(-1)^{l+1}l!}{(\alpha(z)-j+n)^{l+1}}\int_{S_x^*U}
\left(\dd_z^k\s(z)\right)_{\alpha(z)-j,\,l}(x,\xi)\, \dbar_S\xi\right) \nonumber \\
& = &  \hskip 2mm \sum_{l=0}^k \frac{(-1)^{l}\,(l+1)! \
\a^\prime(z)}{(\alpha(z)-j+n)^{l+2}}\int_{S_x^*U}
\ \s^{(k)}(z)_{\a(z)-j,\,l}(x,\xi)\, \dbar_S\xi \nonumber \\
& &  \hskip 10mm  +  \hskip 2mm \sum_{l=0}^k
\frac{(-1)^{l+1}\,l!}{(\alpha(z)-j+n)^{l+1}}\int_{S_x^*U}
\dd_z\left(\s^{(k)}(z)_{\alpha(z)-j,\,l}(x,\xi)\right)\,
\dbar_S\xi, \label{e:ddz}
\end{eqnarray}
where for the second equality we use the property that both of the
factors in each summand on the right-side of \eqref{e:transition}
are holomorphic near $z_0$, and in the notation of
\eqref{componentsderivsigma3}
\begin{equation*}
\left(\dd_z^k\s(z)\right)_{\alpha(z)-j}(x,\xi) = \sum_{r=0}^{k}
\s^{(k)}(z)_{\a(z)-j,\,r}(x,\xi)\log^r[\xi].
\end{equation*}
In that notation the right-side of \eqref{e:transition} for $k$
replaced by $k+1$ reads
\begin{equation}\label{e:transitioncomponents}
\sum_{l=0}^{k+1}
\frac{(-1)^{l+1}l!}{(\alpha(z)-j+n)^{l+1}}\int_{S_x^*U}
\s^{(k+1)}(z)_{\alpha(z)-j,\,l}(x,\xi)\, \dbar_S\xi,
\end{equation}
while on the (co-)sphere $S_x^*U$ where $|\xi|=1$ the identities
of \lemref{inductivederivs} become
\begin{eqnarray*}
\sigma^{(k+1)}_{\alpha(z)-j,\,k+1}(z)(x, \xi) &= & \a^\prime(z)\,
\sigma^{(k)}_{\alpha(z)-j,\,k}(z)(x, \xi),\\[2mm]
\sigma^{(k+1)}_{\alpha(z)-j,\,l}(z)(x, \xi) &= & \a^\prime(z)\,
\sigma^{(k)}_{\alpha(z)-j,\,l-1}(z)(x, \xi)   +
\dd_z(\sigma^{(k)}_{ \alpha(z) -j,\,l} (z)(x, \xi)), \
\ \ \ 1\leq l\leq k,\\[2mm]
\sigma^{(k+1)}_{\alpha(z)-j,\,0}(z)(x, \xi) & = &
\dd_z(\sigma^{(k)}_{ \alpha(z) -j,0} (z)(x, \,\xi)).
\end{eqnarray*}
Substitution of these identities in \eqref{e:transitioncomponents}
immediately shows \eqref{e:transitioncomponents} to be equal to
\eqref{e:ddz}. This completes the proof of
\lemref{lem:transition}.
\end{proof}

\vskip 2mm

Returning to the proof of \thmref{thm:KV2}, from \eqref{ak} and
\eqref{e:cutoffsigmaderivative} and \lemref{lem:transition}  we
now have
\begin{equation*}
a_k=\cutoffint_{T_x^*U} \s^{(k)}(z_0)(x,\xi) \, \dbar\xi,
\end{equation*}
and so the first two lines of \eqref{e:cutoffsigmaz} may be
replaced by $\sum_{k=0}^ K \cutoffint_{T_x^*U}
\s^{(k)}(z_0)(x,\xi)\,\dbar\xi \ \frac{(z-z_0)^k}{k!}+
o((z-z_0)^K)$. Hence \eqref{e:cutoffsigmaz} becomes
\begin{eqnarray}
\cutoffint_{T_x^*U} \sigma(z)(x, \xi)\dbar\xi &  = &
\sum_{k=0}^{K}\cutoffint_{T_x^*U} \s^{(k)}(z_0)(x,\xi) \, \dbar\xi
\
\frac{(z-z_0)^k}{k!}  \ + \ o((z-z_0)^K)\nonumber\\
& & \hskip 10mm  - \ \ {1\over \alpha(z)-\alpha(z_0)}
\int_{S_x^*U}\sigma(z)_{-n}(x, \xi)\, \dbar_S\xi.\label{e:resterm}
\end{eqnarray}

To expand the sphere integral term in \eqref{e:resterm}, since
$\a$ is holomorphic we have in a neighbourhood of each $z_0\in P$
a Taylor expansion $$\alpha(z)-\alpha(z_0)= \sum_{l=1}^L
\frac{\alpha^{(l)}(z_0)}{l!} \, (z-z_0)^l +o(z-z_0)^L$$ and hence
since $\a^\prime(z_0)\neq 0$ an expansion
$$\frac{1}{\alpha(z)-\alpha(z_0)}  =
\frac{1}{\alpha^\prime(z_0)\,(z-z_0)}\,.\, \frac{1}{1 +
\sum_{l=1}^L \frac{\alpha^{(l+1)}(z_0)}{\alpha^\prime(z_0)} \,
\frac{(z-z_0)^l}{(l+1)!} +o(z-z_0)^L}$$
\begin{equation}\label{e:alpha}
 =  \frac{1}{\alpha^\prime(z_0)}\,.\, \frac{1}{(z-z_0)} \ - \
\frac{\alpha^{\prime\prime}(z_0)}{2\,\alpha^\prime(z_0)^2}
 \ + \ \sum_{j=1}^J \b_j(z_0)\, (z-z_0)^j \ + \ o(z-z_0)^J
\end{equation}
with $\b_j(z_0)$ an explicitly computable rational function in
$\a^{(k)}(z_0), \, 1\leq k\leq j+1$ with denominator an integer
power of $\a^{\prime}(z_0)$. On the other hand, since $\alpha(z_0)
-j_0 = -n$, the expansion \eqref{e:sigmazexpansion0} for $j=j_0$
becomes
\begin{equation}\label{e:sigmazexpansion0a}
\sigma(z)_{ \alpha(z) -j_0}(x, \xi) =
\sum_{k=0}^{\o}\sum_{l=0}^k(\s^{(k)}(z_0))_{-n,l}(x,
\xi)\log^l[\xi]\,\frac{(z-z_0)^k}{k!}.
\end{equation}
Since $(\s^{(k)}(z_0))_{-n,l}(x, \xi)\log^l|\xi| =0$ for $l\geq 1$
on $S_x^*U$, we find from the expansions \eqref{e:alpha} and
\eqref{e:sigmazexpansion0a}
$${1\over \alpha(z)-\alpha(z_0)} \int_{S_x^*U}\sigma(z)_{\a(z)-j_0}(x, \xi)\,
\dbar_S\xi=\frac{1}{\alpha^\prime(z_0)}\,.\, \int_{S_x^*U}
\left(\sigma^\prime(z_0)\right)_{-n,0}(x,\xi) \,\dbar_S\xi \,
\frac{1}{(z-z_0)}$$ $$ - \  \sum_{k=0}^K \int_{S_x^*U} {\mathcal
L}_k \left(\sigma(z_0), \sigma^\prime(z_0), \cdots,
\sigma^{(k+1)}(z_0)\right)_{-n,0}(x, \xi)\,
\dbar_S\xi\,\frac{(z-z_0)^k}{k!}
$$
\begin{equation}\label{pole}
+ o\left((z-z_0)^K\right),
\end{equation}
\noi where ${\mathcal L}_k\left(\sigma(z_0), \sigma^\prime(z_0),
\cdots, \sigma^{(k+1)}(z_0)\right)$ is readily seen to have the
form in \eqref{e:L}. In particular, the explicit formulae given
for the first two terms in \eqref{e:alpha} lead to the formula

$${\mathcal L}_0 \left(\sigma(z_0), \sigma^{\prime}(z_0)\right)(x, \xi)
\ =
$$ $$\ \frac{1} {\alpha^\prime(z_0)}\ \int_{S_x^*U}
\sigma^\prime(z_0)_{-n,0} (x,\xi) \ \dbar_S\xi \  - \
\frac{\a^{\prime\prime}(z_0)}{2\,\a^\prime(z_0)^2}\ \int_{S_x^*U}
\sigma(z_0)_{-n} (x,\xi) \ \dbar_S\xi$$

\vskip 2mm

\noi which with the contribution from the $k=0$ finite-part
integral on the right-side of \eqref{e:resterm} gives the stated
constant term in the expansion \eqref{eq:KV2symb}. The next term
up, for example, is

$${\mathcal L}_1 \left(\sigma(z_0), \sigma^\prime(z_0),
\sigma^{\prime\prime}(z_0)\right)(x, \xi) \ =
$$ $$\ \frac{1} {\alpha^\prime(z_0)}\ \int_{S_x^*U}
\sigma^{\prime\prime}(z_0)_{-n,0} (x,\xi) \ \dbar_S\xi \  - \
\frac{\a^{\prime\prime}(z_0)}{2\,\a^\prime(z_0)^2}\ \int_{S_x^*U}
\sigma^{\prime}(z_0)_{-n,0} (x,\xi) \ \dbar_S\xi$$
$$ + \ \frac{3\a^{\prime\prime}(z_0)^2 -
2\a^{\prime\prime\prime}(z_0)\a^\prime(z_0)}{12\,\a^\prime(z_0)^3}\
\int_{S_x^*U} \sigma(z_0)_{-n} (x,\xi) \ \dbar_S\xi.$$

\vskip 2mm

When $\alpha(z) = qz+b$ with $q\neq 0$ the right-side of
\eqref{e:alpha} is $ \frac{1}{q \, (z-z_0)}$ and so from
\eqref{e:sigmazexpansion0} one then has ${\mathcal
L}_k\left(\sigma(z_0), \sigma^\prime(z_0), \cdots,
\sigma^{(k+1)}(z_0)\right)= \frac{\sigma^{(k+1)}(z_0)}{
q\,(k+1)!}$ and so \eqref{eq:KV2symbq} follows.

If $z_0\notin P$ then $\a(z)\in\C\backslash\Z$ and so the
log-polyhomogeneous symbols ${\mathcal L}_k$ in \eqref{pole} then
have non-integer order and hence have no component of degree $-n$,
and therefore vanish. Likewise the pole in\eqref{pole} vanishes
and so \eqref{eq:KV2symb} simplifies,  in this case, to
\eqref{eq:KV2symbhol}. Alternatively, this can be seen in a
simpler more direct way by using the linearity of the finite-part
integral over log-polyhomogeneous symbols of non-integer order
applied to the Taylor expansion of the symbol at $z_0$. (Indeed,
in this case the term $j=j_0$ in \eqref{e:cutoffsigmaz} does not
need to be treated separately from the sum in the previous line
and \eqref{ak} holds by linearity, from which
\lemref{lem:transition} may then be inferred and now including the
case $j=j_0$.)

This shows the first part of the theorem.

\vskip 2mm

For the second part we use a partition of unity $\{(U_i,\phi_i) \
| \ i\in J\}$  such that for $i,j\in J$ there is an $l_{ij}\in J$
with $\supp(\phi_i) \cup \supp(\phi_j) \subset U_{ij} :=
U_{l_{ij}}$. We suppose trivialisations of $\pi: E\to M$ over each
open set $U_i$. Then, with $U_{ij}$ identified with an open subset
of $\R^n$, one has $A(z) = \sum_{i,j}\phi_i A(z) \phi_j$ where
$\phi_i A(z) \phi_j = {\rm Op}(\s_{(ij)}(z))$ is the localization
of $A$ over $U_{ij}$ with amplitude
$$\s_{(ij)}(z)(x,y,\xi)\in \CS^{\a(z)}(U_{ij}\times U_{ij}, V)$$ a
local holomorphic family of symbols in $(x,y)$ form. Each
finite-part integral $\cutoffint_{T_x
U_{ij}}\s_{(ij)}(z)(x,x,\xi)\ \dbar\xi $ is well defined outside
$P$, since $A(z)$ has non-integer order for those values of $z$.
Using the linearity there of the canonical trace functional it
follows that for $z\notin P$
\begin{equation*}
\TR(A(z))  = \sum_{i,j}\int_{U_{ij}} \cutoffint_{T_x
U_{ij}}\tr\left(\s_{(ij)}(z)(x,x,\xi)\right) \ \dbar\xi\,dx,
\end{equation*}
where $\tr$ is the trace on $\End(V)$, allowing \eqref{eq:KV2symb}
to be applied to each of the summands defined over the
trivialising charts. Each locally defined coefficient in the
Laurent expansion is seen by holomorphic continuation to define a
global density on $M$ in the way explained in
\propref{prop:densities}. The first part of the theorem therefore
yields that ${\rm TR}(A(z))$ is meromorphic with simple poles in
$P$ and since
\begin{equation}\label{e:sigmak}
\sigma^{(k)}_{A(z_0)} = \sigma_{A^{(k)}(z_0)}
\end{equation}
the identity \eqref{eq:KV2op} now follows from the  formula
\eqref{eq:KV2symb} applied to each localization.

\vskip 1mm

The fact that the coefficients of $\frac{(z-z_0)^k}{k! }$ in the
Laurent expansions of the meromorphic maps $z\mapsto
\cutoffint_{T_x^ *U}\sigma(z)$ and  $z\mapsto {\rm TR}(A(z))$
correspond to the finite part at $z=z_0$ of their derivative at
order $k$ follows from the general property for a meromorphic
function $f$ on an open set $W\subset \C$ with Laurent expansion
around $z_0$ given by $f(z)=\sum_{j=1}^J
\frac{b_j}{(z-z_0)^j}+\sum_{k=0}^K a_k
\frac{(z-z_0)^k}{k!}+o((z-z_0)^K)$ that
\begin{equation}\label{e:mero}
{\rm fp}_{z=z_0} f^{(k)}(z)= a_k.
\end{equation}
Combined with the equality $\dd^k_z \TR(A(z)) = \TR(A^{(k)}(z))$
valid for $z\not\in P$ we reach the conclusion.

Since the formulas \eqref{linearorder}, \eqref{eq:KV2ophol} now
follow from \eqref{eq:KV2symbq} and \eqref{eq:KV2symbhol}, this
ends the proof of the theorem.
\end{proof}

\vskip 5mm

In passing  from the local formula \eqref{eq:KV2symb} to the
global formula \eqref{eq:KV2op} in the  proof of \thmref{thm:KV2}
we have implicitly used the following fact, yielding the Laurent
coefficients to be global densities on $M$ which can be
integrated.
\begin{prop}\label{prop:densities}
Let $c_k(x)$ denote the coefficient of $\frac{(z-z_0)^k}{k!}$ in
the Laurent expansion \eqref{eq:KV2symb}. Then $c_k(x)\,dx$ is
defined independently of the choice of local coordinates on $M$.
\end{prop}
\begin{proof}
By formula \eqref{e:mero}, the coefficient $c_k(x)$ of
$\frac{(z-z_0)^k}{k!}$ in the Laurent expansion \eqref{eq:KV2symb}
with $\sigma(z)(x,\cdot)=\sigma_{A(z)}(x,\cdot)$  is identified
with the finite part at $z_0$ of the $k$-th derivative of the  map
$$ z\mapsto I_{A(z)}(x):=\cutoffint_{T_x^* U} \sigma_{A(z)}(x,
\xi)\, \dbar\xi$$ i.e. $c_k(x)= {\rm fp}_{z=z_0}
I_{A^{(k)}(z)}(x).$ For $z\notin P$ the property
\eqref{e:TRdensity} holds for the finite part integral
$I_{A(z)}(x)$ as well as  for the finite part integrals
$I_{A^{(k)}(z)}(x)$ since the order of $A^{(k)}$ differs from that
of $A(z)$ by an integer.

The map $z\mapsto I_{A^{(k)}(z)}(x)$ has a Laurent expansion
$I_{A^{(k)}(z)}(x)=\sum_{j=1}^{k+1} \frac{b_j(x)}{(z-z_0)^j}+
\sum_{k=0}^Kc_k(x)\,\frac{(z-z_0)^k}{k!} +o((z-z_0)^k)$ and
$(z-z_0)^{k+1} I_{A^{(k)}(z)}(x)$ can be extended to a holomorphic
function in a small ball centered at  $z_0$ with value
$b_{k+1}(x)$ at $z_0$. Since  property \eqref{e:TRdensity} holds
for $ I_{A^{(k)}(z) }(x)$ outside $z_0$  in this ball, it holds
for the holomorphic extension on the whole ball and hence for
$b_{k+1}(x)$. Using \eqref{e:TRdensity}, we deduce that
$b_{k+1}(x)\, dx$ is defined independently of the choice of local
coordinates on $M$ and  so is the difference
$\left(I_{A^{(k)}(z)}(x)-
\frac{b_{k+1}(x)}{(z-z_0)^{k+1}}\right)\, dx$  for any $z$ outside
$z_0$ in a small ball centered at  $z_0$.  Iterating    this
argument, one shows recursively on the integer $1\leq J\leq k$
that  $\left(I_{A^{(k)}(z)}(x)-\sum_{j=1}^{k+1-J}
\frac{b_j(x)}{(z-z_0)^j}\right)\, dx$ is defined independently of
the choice of local coordinates on $M$ in a  small ball centered
at $z_0$. Consequently,  the finite part $\left({\rm
fp}_{z=z_0}I_{A^{(k)}(z)}(x)\right)\, dx$ at  $z_0$ is also
defined independently of the choice of local coordinates. Since
this finite part coincides with $k!\, c_k(x)$, we have  that
$c_k(x)\,dx$ is defined independently of the choice of local
coordinates on $M$.
\end{proof}

\vskip 2mm

Examining the singular and constant terms in the expansions of
\thmref{thm:KV2} we have the following corollaries.

\vskip 2mm

 \noi First, the singular term yields the known identification of
the residue trace with complex residue of the canonical trace,
derived in \cite{Gu}, \cite{Wo}, \cite{KV}. With the assumptions
of \thmref{thm:KV2}:

\begin{cor}\label{cor:KV}
The map $z\mapsto \cutoffint_{T_x^*U} \sigma(z)( x,\xi)
\,\dbar\xi$ is meromorphic with at most a simple pole at $z_0\in
P$ with complex residue

\begin{equation} \label{eq:KVsymb}
{\rm Res}_{z=z_0}\cutoffint_{T_x^*U} \,\sigma(z) (x, \xi) \,
\dbar\xi = -{1\over \alpha^\prime( z_0 )}\ \int_{S_x^*U}
\sigma(z_0)_{-n} (x,\xi) \ \dbar_S\xi.
\end{equation}

\vskip 3mm

For the holomorphic family $z\mtoo A(z)$ of $\pdo$s parametrized
by $W$, the form $\ \frac{1} {\alpha^\prime(z_0)}\ \int_{S_x^*U}
\left(\sigma_{A(z_0)} \right)_{-n} (x,\xi) \ \dbar_S\xi\, dx$
defines a global density on the manifold $M$ and the  map
$z\mapsto {\rm TR}(A(z)):= \int_M \, dx\, \TRx(A(z))$ is a
meromorphic function with at most a simple pole at $z_0\in P$ with
complex residue

\begin{equation} \label{eq:KVop}
{\rm Res}_{z=z_0}\TR(A(z) )= -{1\over \alpha^\prime( z_0 )}\,
\res\left( A(z_0) \right).
\end{equation}
\end{cor}

\vskip 2mm

Thus, consequent to \propref{prop:densities}, one infers here the
global existence of the residue density for integer order
operators from  the existence of the canonical trace density for
non-integer order operators and holomorphicity.

\vskip 5mm

\noi On the other hand, the constant term provides a `defect
formula' for finite part integrals.

\vskip 2mm

With the assumptions of \thmref{thm:KV2}:

\begin{thm}\label{cor:defect}
For  a holomorphic family of symbols $z\mapsto \sigma(z)\in \CS(U,
V)$ parametrized by a domain $W\subset \C$ and for any  $x\in U$,

\begin{eqnarray}
{\rm fp}_{z=z_0}\cutoffint_{T_x^*U}  \sigma(z) (x, \xi)\, \dbar\xi
\ &=& \ \cutoffint_{T_x^*U}\sigma(z_0) (x, \xi)\,\dbar\xi \ - \
{1\over \alpha^\prime( z_0 )} \, \int_{S_x^*U}
\sigma^\prime(z_0)_{-n,0} (x,\xi) \ \dbar_S\xi.
 \nonumber\\
&&\hskip 10mm  + \ \
\frac{\a^{\prime\prime}(z_0)}{2\,\a^\prime(z_0)^2}\ \int_{S_x^*U}
\sigma(z_0)_{-n} (x,\xi) \ \dbar_S\xi. \label{eq:defectsymb}
\end{eqnarray}

\vskip 1mm

\noi For the holomorphic family $z\mapsto A(z)\in \Cl(M, E)$ of
$\pdo$s parametrized by $W\subset \C$,

\vskip 1mm

\begin{eqnarray}
{\rm fp}_{z=z_0}\TR(A(z) ) \ &=& \ \int_M dx\ \left(\TRx(A(z_0)) \
- \ \frac{1}{\alpha^\prime(z_0)}\, \resxo(A^\prime(z_0))\,\right)
\nonumber \\
& &\hskip 10mm + \ \
\frac{\a^{\prime\prime}(z_0)}{2\,\a^\prime(z_0)^2}\
\res\,\left(A(z_0) \right) \label{eq:defectop}
\end{eqnarray}
\end{thm}

\vskip 1mm

\begin{rk}
Since $\a$ is non-critical on $P$, from \propref{prop:diffsymb} if
$z_0\in P$ the operator $A^\prime(z_0)\in \Cl^{\a(z_0),1}(M,E)$ in
equation \eqref{eq:defectop} is not classical.
\end{rk}

\begin{rk}
If ${\rm res}_{x,0}(A(z_0))=0$ then ${\rm fp}_{z=z_0}\TRx(A(z) )
=\lim_{z\to z_0} \TRx(A(z))$. If this holds for all $x\in M$, then
$\TR(A(z))$ is holomorphic at $z_0$ and ${\rm fp}_{z=z_0}\TR(A(z)
) =\lim_{z\to z_0} \TR(A(z))$.
\end{rk}

\vskip 3mm

\noi One therefore has the following statement on the existence of
densities associated to the local canonical trace.

\begin{thm}\label{t:gendensity}
With the assumptions of \thmref{thm:KV2}, for a holomorphic family
$z\mapsto A(z)\in \Cl(M, E)$ parametrized by a domain $W\subset
\C$, and irrespective of the order $\a(z_0)\in\Rf$ of $A(z_0)$
\begin{equation}\label{anyorderdensity}
\left(\TRx(A(z_0)) \ - \ \frac{1}{\alpha^\prime(z_0)}\,
\resxo(A^\prime(z_0))\,\right)dx
\end{equation}defines a global density on $M$
which integrates on $M$ to ${\rm fp}_{z=z_0} {\rm TR}(A(z))$. If
$\a(z_0)\notin\Z$ then \eqref{anyorderdensity} reduces to the
canonical trace density on non-integer order classical $\pdo$s of
\cite{KV}.
\end{thm}

\vskip 2mm \noi Though this follows on the general grounds of
\propref{prop:densities}, we have, for completeness, given a
direct proof of \thmref{t:gendensity} in Appendix A. This
specializes to give the previously known existence of the
canonical trace on non-integer order $\pdo$s, recalled in Section
(1.2).

\vskip 2mm

With the assumptions of \thmref{thm:KV2}:

\begin{thm}\label{prop:globaldefect}
Let $z\mapsto A(z)\in \Cl(M, E)$ be a holomorphic   family of
classical $\pdo$s parametrised by $W\subset \C$ and let $z_0\in
W$. If either
$$\TRx(A(z_0))\, dx=\left( \cutoffint_{T_x^*M} {\rm tr}_{x}
(\sigma_{A(z_0)}) (x, \xi)\, d\xi\right)\, dx$$ or
$${\rm res}_x(A^\prime(z_0))\, dx:= \int_{S^*_xM}\,  {\rm tr}_{x}
\left(\left(\sigma_{A^\prime(z_0)}\right)_{-n} (x, \xi)
\right)\,d_S\xi\, dx$$ defines a global density  on $M$, then
${\rm TR}\left(A(z_0)\right)$ and ${\rm
res}\left(A^\prime(z_0)\right)= \int_M {\rm
res}_{x,0}(A^\prime(z_0))\, dx$ are both well defined. The
following defect formula then holds
\begin{equation}\label{e:globaldefect}
{\rm fp}_{z=z_0}\TR(A(z) ) \ = \ {\rm TR}(A(z_0)) \ -\
\frac{1}{\alpha^\prime(z_0)}\,{\rm res}\left(A^\prime(z_0)\right).
\end{equation}
\vskip 2mm \noi This holds in the following cases:

\vskip 2mm

\noi {\rm (i)} \ If $A(z_0)\in \Cl^{\a(z_0),\,0}(M,E)$ satisfies
one of the cases {\rm (1), (2)} or {\rm (3)} of
\propref{p:TRdefined} then ${\rm TR}(A(z_0))$ is defined and
\eqref{e:globaldefect} holds. In case (1) this reduces to
\begin{equation}\label{e:resAprimezero1}
{\rm fp}_{z=z_0}\TR(A(z) )= {\rm TR}(A(z_0)).
\end{equation}

 \vskip 2mm \noi {\rm
(ii)} \ If ${\rm res}_{x,0}\left(A^\prime(z_0)\right)=0$ for all
$x\in M$ then $\TR(A(z))$ is holomorphic at $z_0\in W$, so that
${\rm fp}_{z=z_0}\TR(A(z) ) = \lim_{z\to z_0} {\rm TR}(A(z))$, and
\eqref{e:resAprimezero1} holds.

\vskip 2mm \noi {\rm (iii)} \ If $A(z_0)$ is a differential
operator, and more generally whenever $\TRx(A(z_0))=0$ for all
$x\in M$, \eqref{e:globaldefect} reduces to
\begin{equation}\label{e:TRAzero1}
{\rm fp}_{z=z_0}\TR(A(z) )=
-\frac{1}{\alpha^\prime(z_0)}\res\left(A^\prime(z_0)\right).
\end{equation}
\end{thm}
\begin{rk} \eqref{e:globaldefect} can hold with both summands on the
right-side of the equation non-zero. See  \exref{rk:oddorder}.
\end{rk}
\begin{proof}
The first statement is consequent to \thmref{t:gendensity}. Since
$\TRx(A(z_0))\, dx$ then defines a global density the
transformation rule for finite part integrals in
\propref{prop:noninv} implies that
$$\int_{S_x^*M} \, {\rm tr}_x\left(
\sigma_{A(z_0)}\right)_{-n,\,0} (x, \xi) \, \log \vert
C^{-1}\xi\vert \, \dbar_S\xi = 0\quad \forall \, C\in GL_n( \C)$$
and hence (taking $C=\l\cdot  I$, $\la\in\C$) that $${\rm
res}_{x,0}(A(z_0)) := \int_{S_x^*M} \dbar_S\xi\, {\rm tr}_x\left(
\sigma_{A(z_0)}\right)_{-n,\,0} (x, \xi)=0.$$ Equation
\eqref{e:globaldefect} now follows from \eqref{eq:defectop}. Parts
(i), (ii), (iii) are now obvious in view of \propref{p:TRdefined}
and \propref{prop:cutoffdiff} and the vanishing of the residue
trace on non-integer order operators and on differential
operators.
\end{proof}

\vskip 5mm

\section{Application to the Complex Powers}

An operator $Q\in \Ell(M, E)$ of positive order  is called {\it
admissible } if there is a proper subsector of $\Cf$ with vertex 0
which contains the spectrum of the leading symbol $\sigma_L(Q)$ of
$Q$. Then there is a half line $L_\theta=\{re^{i\theta}, r>0\}$ (a
spectral cut) with vertex $0$ and  determined by an Agmon angle
$\theta$ which does not intersect the spectrum of $Q$. Let
$\Ell_{\ord>0}^{ adm}  (M, E)$  denote the subset of admissible
operators in $\Ell (M, E)$ with positive order.

Let $Q\in \Ell^{adm}_{\ord >0}(M, E)$ with spectral cut
$L_\theta$. For Re $z < 0$, the complex power $Q_\theta^z$ of $Q$
is a bounded operator on any space $H^s(M, E)$ of  sections of $E$
of Sobolev class $H^s$  defined by the contour integral:
\begin{equation}\label{eq:complexpower}
Q_\theta^z={i\over 2\pi} \int_{C_{\theta}} \l^z (Q-\l I)^{-1} d\l
\end{equation}
where  $C_{\theta}= C_{1, \theta,r}\cup C_{2, \theta,r}\cup
C_{3,\theta,r}$. Here  $r$ is a sufficiently small positive number
and $C_{1,\theta,r}= \{\l=  \vert \lambda \vert e^{i\theta} \ | \
+\infty> \vert \lambda \vert \geq r\}$, $C_{2,\theta,r}= \{\l= r
e^{i \phi } \ | \  \theta\geq\phi\geq \theta-2\pi \}$ and
$C_{3,\theta,r}= \{\l=  \vert \lambda\vert e^{ i (\theta-2\pi)  }
\ | \ r \leq \vert \lambda\vert<+\infty\}$.
%%$C_{4,\theta,r}=
%%\{\l= R e^{i \theta}, \theta\geq\phi\geq
%%-2\pi+\theta \}$
Here $\l^z= \exp (z\log \l)$ where $\log \l= \log \vert \l\vert +
i \theta$ on $C_{1, \theta,r}$ and $\log \l= \log \vert \l\vert  +
i (\theta-2\pi)$ on
$C_{3, \theta,r}$. \\ \\
For $k\in \N$ the complex power $Q^z$ is then extended to the half
plane Re $z<k$ via the relation \cite{Se1} $$ Q^k Q_\theta^{z-k}=
Q_\theta^z .$$ The definition  of a complex power  depends in
general on the choice of $\theta$ and yields for any $z\in \C$  an
elliptic operator $Q_\theta^z$ of order $z\cdot \ord(Q)$. In spite
of this $\theta$-dependence, we may occasionally omit it in order
to simplify notations.
\begin{rk} For $z=0$
$$Q_\theta^0= I - \Pi_Q$$
where $\Pi_Q$ is the smoothing operator projection
$$\Pi_Q =
\frac{i}{2\pi} \int_{C_0} (Q-\lambda I)^{-1} \, d\lambda$$ with
$C_0$ a contour containing the origin but no other element of
${\rm spec}(Q)$, with range the generalized kernel $\{\psi\in
\Ci(M, E)\ | \ Q^N\psi=0 \ \ {\rm for} \ {\rm some} \ N\in \N\}$
of $Q$. (See \cite{Bu}, \cite{Wo}, presented recently in
\cite{Po}).
\end{rk}
Let   $Q\in \Ell_{\ord >0}^{ adm} (M, E)$ be of  order $q$ with
spectral cut $L_\theta$. For arbitrary $k\in \Z$, the map $z\to
Q^z_\theta$  defines a holomorphic function from $\{z\in \C, {\rm
Re}z<k\}$ to the space ${\mathcal L} \left(H^s(M, E)\to
H^{s-k\cdot q }(M, E)\right)$ of bounded linear maps and we can
set
$$\log_\theta Q:= \left[{\partial \over
\partial z} Q_\theta^z\right]_{z=0}. $$
From \eqref{e:derivsigma}, in a local trivialisation
$E_{\vert_U}\simeq U\times V$ of $E$ over an open set $U$ of $M$
the  symbol of $\log_\theta Q$ reads $\sigma_{\log_\theta Q}(x,
\xi)= \ord(Q) \log \vert \xi\vert {\rm Id} + \rho(x, \xi)$ with
$\rho \in \Cl^0(U,V)$, and so $\log_\theta Q \in \Cl^{0,1}(M,E)$
has order zero and log degree one. The logarithmic dependence is
slight, for $P,Q\in \Ell_{\ord>0}^{ adm} (M, E)$, of non zero
order $p, q$ respectively and admitting spectral cuts $L_\theta$
and $L_\phi$ we have ${\log_\theta P\over p }- {\log_\phi Q\over q
}\in \Cl^0(M, E). $ More generally, higher derivatives of the
complex powers have symbols with polynomial powers of $\log \vert
\xi\vert$ and it follows from \propref{prop:diffsymb} that
\begin{equation}\label{e:logk}
\log^k_\theta Q:= \left[{\partial^k \over
\partial z^k} Q_\theta^z\right]_{z=0} \in \Cl^{0,k}(M,E).
\end{equation}

\vskip 1mm

Theorem \ref{thm:KV2} leads to the following Laurent expansion.

\begin{thm}\label{thm:KV3}
Let $Q\in \Ell_{ord >0}^{adm}(M, E)$ with spectral cut $\theta$
and of order $q$ and let $A\in \Cl^{\a}(M, E)$. On the half plane
${\rm Re}(z)>\frac{\alpha+n}{q}$ the local Schwartz kernel $K_{A\,
Q_\theta^{-z}}(x,y)$ of $AQ_\theta^{-z}$ is well defined and
holomorphic and the restriction to the diagonal $ K_{A\,
Q_\theta^{-z}}(x,x)\,dx = \int_{T_x^*M} \sigma_{A\,
Q_\theta^{-z}}(x, \xi) \, \dbar\xi\,dx$ defines a global density,
an element of $\Ci(M, {\rm End}(E))$. There is a meromorphic
extension of $K_{A\, Q_\theta^{-z}}(x,y)$ to all $z\in\C$
\begin{equation}\label{eq:Kmero}
K_{A\, Q_\theta^{-z}}(x,x)\vert^{{\rm mer}}:= \cutoffint_{T_x^*M}
\sigma_{A\, Q_\theta^{-z}}(x, \xi) \, \dbar\xi
\end{equation}
with at most simple poles, each of which is  located in
$P:=\{\frac{\alpha-j}{q}\,|\, j\in [-n, \infty[ \ \cap\ \Z\}$. For
any $x\in M$, we have for $z$ near $\frac{\alpha-j}{q}\in P$
\begin{eqnarray} \label{eq:KV3symb}
 \  \  && K_{A\, Q_\theta^{-z}}(x,x)\vert^{{\rm mer}} \ \ = \ \
\frac{1}{ q }\, \int_{S^*_xM}\,  \left(\sigma_{A\,
Q_\theta^{(j-\alpha)/q}}\right)_{-n} (x, \xi) \,\dbar_S\xi \cdot
\frac{1}{
(z-\frac{\alpha-j}{q})}\nonumber\\[1mm] & + &  \ \sum_{k=0}^K\frac{(-1)^{k}}{k!}
\left(z-\frac{\alpha-j}{q}\right)^k\nonumber\\
& &  \times \, \left( \cutoffint_{T_x^*U}\, \sigma_{A\,
Q_\theta^{(j-\alpha)/q}\log_\theta^k Q}(x, \xi)\,\dbar\xi
%\TRx(A\,Q^{-\frac{\alpha-j}{q}}\, \log_\theta^{(k)}Q  )
-\frac{1}{q(k+1)}\, \int_{S^*_xM}\,  \left(\sigma_{A\,
Q_\theta^{(j-\alpha)/q}\,\log_\theta^{k+1}Q}\right)_{-n,0} (x,
\xi) \,\dbar_S\xi\right) \nonumber\\[3mm]  & + &  \ \
 o\left(\left(z-\frac{\alpha-j}{q}\right)^K\right).
\end{eqnarray}

\noi It follows that the map $z\mapsto {\rm TR}(A\,
Q_\theta^{-z}):= \int_M \,{\rm tr}_x\left(K_{A\,
Q_\theta^{-z}}(x)\vert^{{\rm mer}} \right) $ is a meromorphic
function with no more than simple poles located in $ P$, and for
$z$ near $\frac{\alpha-j}{q}\in P$
\begin{eqnarray}
& & \TR(A\, Q_\theta^{-z} )  \ \  = \ \ \frac{1}{q}\, {\rm
res}(A\, Q_\theta^{\frac{j-\alpha}{q}})\cdot \frac{1}{
(z-\frac{\alpha-j}{q})}\nonumber \\[2mm]
& & + \ \ \sum_{k=0}^K\frac{(-1)^{k}}{k!}
\left(z-\frac{\alpha-j}{q}\right)^k\nonumber\\
& & \hskip 5mm \times \,\int_M dx \left(
\TRx(A\,Q_\theta^{\frac{j-\alpha}{q}}\log_\theta^k Q) -
\frac{1}{q(k+1)} {\rm
res}_{x,0}(A\,Q_{\th}^{\frac{j-\alpha}{q}}\,\log_\theta^{k+1}Q)\right)
\nonumber\\[3mm]
 & & + \ \
  o\left(\left(z-\frac{\alpha-j}{q}\right)^K\right).\label{eq:KV3op}
\end{eqnarray}
If $z_0\notin P$ then $\TR(A\, Q_\theta^{-z} )$ is holomorphic at
$z_0$ and for $z$ in a small enough neighbourhood of $z_0$
\begin{equation} \label{eq:KV3op2}
\TR(A\, Q_\theta^{-z} )  \   =  \ \sum_{k=0}^K\frac{(-1)^{k}}{k!}
 \ \TR(A\,Q_\theta^{z_0}\log_\theta^k
Q)\,\frac{(z-z_0)^k}{k!} \ + \
  o\left(\left(z-z_0\right)^K\right).
\end{equation}

\end{thm}

\begin{proof}
Since $\sigma(z):= \sigma_{A\, Q_{\theta}^{-z}}$ has order
$\alpha(z)= \alpha-q\, z$, \eqref{eq:KV2symbq} of Theorem
\ref{thm:KV2} can be applied to equation \eqref{eq:Kmero}. Using
\eqref{e:logk} and \exref{example}, this yields
(\ref{eq:KV3symb}). Applying the fibrewise trace ${\rm tr}_x$ and
integrating over $M$ yields equation (\ref{eq:KV3op}). Equation
\eqref{eq:KV3op2}, to which \eqref{eq:KV3op} reduces when
$\a\notin\Z$, as the operators inside the local residue traces
then have non-integer order, follows from \eqref{eq:KV2ophol};
that $ \TRx(A\,Q_\theta^{z_0}\log_\theta^k Q)\, dx$  defines a
global density on $M$ in this case is known from \cite{Le}.
\end{proof}

Because of the identity with the generalized zeta-function
$$\zeta_\theta(A, Q,
z) \ =\  {\rm TR}\left(A \, Q_\theta^{-z}\right)$$ the expansion
\eqref{eq:KV3op} is of particular interest near $z=0$, owing to
the role of the Laurent coefficients there in geometric analysis.
\begin{thm}\label{t:TRatzero}
If $\ord(A) = \a \in [-n, \infty[ \ \cap\ \Zf$ then $0\in P$ and
one then has near $z=0$
\begin{eqnarray}
& &\zeta_\theta(A, Q,
z) \ = \ \frac{1}{q}\, {\rm res}(A)\cdot \frac{1}{z}\nonumber \\[2mm]
& & \hskip 5mm + \ \ \int_M dx\, \left( \TRx(A) - \frac{1}{q}\,
{\rm res}_{x,0}(A\,\log_\theta Q)\right) \ - \ {\rm tr}
(A\,\Pi_Q)\,
\nonumber\\[2mm]
& & \hskip 5mm + \ \ \sum_{k=1}^K (-1)^k\,\frac{z^k}{k!}\,\nonumber\\
& & \hskip 10mm \times \int_M dx\, \left( \TRx(A\log_\theta^k Q) -
\frac{1}{q(k+1)}\, {\rm res}_{x,0}(A\,\log_\theta^{k+1}Q)\right) \
- \ {\rm tr} (A\, \log_\theta ^k Q\,\Pi_Q)\,
\nonumber\\[2mm]
& & \hskip 5mm + \ \
 o (z^K).\label{eq:KV3zerotrace}
\end{eqnarray}
If $\a \notin [-n, \infty[ \ \cap\ \Zf$ then $\zeta_\theta(A, Q,
z)$ is holomorphic at zero and one has for  $z$ near zero
\begin{equation}
\zeta_\theta(A, Q, z) =  \sum_{k=0}^K
(-1)^k\,\left(\TR(A\log_\theta^k Q)  -  {\rm tr} (A\, \log_\theta
^k Q\,\Pi_Q)\right)\,\frac{z^k}{k!}  \ + \ \ o
(z^K).\label{eq:KV3zerotracenoninteger}
\end{equation}
\end{thm}
\begin{rk}
The formula \eqref{eq:KV3zerotracenoninteger} can also be deduced
from exact formulas for the case $\a\notin\Z$ in \cite{Gr1}
Sect(3). All formulas presuppose the existence shown in \cite{KV},
\cite{Le} of the canonical trace for non-integer order $\pdo$s
with log-polyhomogeneous symbol.
\end{rk}
\begin{proof}
The assumption $\a \in [-n, \infty[ \ \cap\ \Zf$ means that $(\a
-j_0)/q =0$ for some $j_0 \in [-n, \infty[ \ \cap\ \Zf$. Hence
\eqref{eq:KV3zerotrace} is almost obvious from \eqref{eq:KV3op};
the subtle point is to take care to replace
$Q_{\th}^{-\frac{\alpha-j}{q}} = Q_{\th}^0$ by $I - \Pi_Q$, see
Remark 2.1. Since the spectral projection $\Pi_Q$ is a smoothing
operator the term  $\TRx(A\, \log_\theta ^k Q\,\Pi_Q)\,dx$ is an
ordinary integral valued density and globally defined, yielding
the term $ {\rm tr} (A\, \log_\theta ^k Q\,\Pi_Q)$. The formula
\eqref{eq:KV3zerotracenoninteger} for $A$ of non-integer order (to
which \eqref{eq:KV3zerotrace} reduces in this case) is immediate
from \eqref{eq:KV3op2}.
\end{proof}

We denote the coefficient of $(z-\frac{\alpha-j}{q})^k/k!$  in the
Laurent expansion of the generalized zeta function at
$\frac{\alpha-j}{q}\in P$ by $ \zeta^{(k)}_\theta(A, Q,
\frac{\alpha-j}{q})$. In the case $k=0$,  we use the simpler
convention of writing the constant term $\zeta^{(0)}_\theta(A, Q,
\frac{\alpha-j}{q}) := {\rm fp}_{z=
\frac{\alpha-j}{q}}\,\zeta_\theta(A, Q, z)$  as $\zeta_\theta(A,
Q, \frac{\alpha-j}{q})$.
 When $A=I$ write $\zeta_\theta( Q,
\frac{\alpha-j}{q}):=\zeta_\theta(I, Q, \frac{\alpha-j}{q}).$

\vskip 3mm

\begin{cor}\label{thm:zetamer}
For any operator $A\in \Cl(M, E)$,
\begin{equation}\label{eq:zetaatzero}
\zeta_\theta(A, Q, 0 ) = \int_M dx \left( \TRx(A) - \frac{1}{q}\,
{\rm res}_{x,0}(A\,\,\log_\theta Q)\right)  \ - \ {\rm tr}
(A\,\Pi_Q).
\end{equation}
More generally, for any non-negative integer $k$
\begin{eqnarray}\label{eq:zetakatzero}
\zeta_\theta^{(k)}(A,Q, 0) \ & = &\ (-1)^k\, \int_M dx\,
\left(\TRx (A\, \log_\theta^k Q) - \frac{1}{q\, (k+1)}\,{\rm
res}_{x,0}(A\,\log_\theta^{k+1} Q )\right) \
\nonumber \\[2mm] && \hskip 20mm + \ \  (-1)^{k+1}\, {\rm tr}\left(A\,
\, \log_\theta^k Q\, \Pi_Q\right) \ .
\end{eqnarray}
If $A$ has integer order $\a \in [-n, \infty)\,\cap\ \Zf$ then  \
\begin{equation}\label{eq:zetaatzj}
\zeta_\theta\left(A, Q,  \frac{\alpha-j}{q} \right)\nonumber\\
= \int_M dx \left( \TRx(A\,Q_\theta^{-\frac{\alpha-j}{q}}) -
\frac{1}{q}\, {\rm
res}_{x,0}(A\,Q_\th^{-\frac{\alpha-j}{q}}\,\log_\theta Q)\right)
\end{equation}
\end{cor}

\vskip 3mm

\noi Applied to the complex powers, the general statement on the
existence of densities associated to the canonical and residue
traces of \thmref{t:gendensity} now states that independently of
the order of $A\in\Cl(M,E)$,
$$\left(\TRx(A)
- \frac{1}{q}\, {\rm res}_{x,0}(A\,\,\log_\theta Q)\,\right)dx$$
always defines a global density on $M$.

\vskip 2mm

\noi If $A$ has non-integer order this reduces to the KV canonical
trace density and (by \eqref{eq:KV3zerotracenoninteger}) the
identity \eqref{eq:zetakatzero} loses its residue defect term and
one then has the known formula (cf. \cite{Gr1} Cor. (3.8))
\begin{equation}
 \zeta^{(k)}_\theta(A, Q, 0) \ = \
(-1)^k\,\TR(A\log_\theta^k Q)\ - (-1)^k\,{\rm tr} (A\, \log_\theta
^k Q\,\Pi_Q).
\end{equation}

\vskip 2mm   Applying \thmref{prop:globaldefect} to the zeta
function at $z=0$ yields the following refinement of
\eqref{eq:zetaatzero}.

\begin{thm}\label{t:globaldefectzeta}
Let $Q\in \Ell_{ord >0}^{adm}(M, E)$ be  a classical $\pdo$ with
spectral cut $\theta$ and of order $q$ and let $A\in
\Cl^{\a}(M,E)$ be  a classical $\pdo$ of order $\a$. If either
$\TRx(A)\, dx$ or ${\rm res}_x(A\,\log_{\th}Q)\, dx$ defines a
global density on $M$, then $\z_{\th}(A,Q,z)$ is holomorphic at
$z=0$, ${\rm TR}(A)$ and ${\rm res}(A\,\log_{\th}Q)$ both exist,
and one has
\begin{equation}\label{e:zetaglobaldefect}
\z_{\th}(A,Q,0) \ = \ {\rm TR}(A) \ - \ \frac{1}{q}\,{\rm
res}(A\,\log_{\th}Q)  \ - \ {\rm tr} (A\,\Pi_Q)\,.
\end{equation}
\end{thm}
\begin{proof}
If $\TRx(A)\, dx$ defines a global density then $\res(A)$
vanishes, as accounted for in the proof of
\thmref{prop:globaldefect}, and so $\z_{\th}(A,Q,z)$ is
holomorphic at $z=0$. The formula is obvious from
\eqref{e:globaldefect}.
\end{proof}

\vskip 2mm

Notice that the assumptions of \thmref{t:globaldefectzeta} also
force $\res(A)=0$.

\vskip 1mm

The situation of \thmref{t:globaldefectzeta} can be seen to hold
for certain combinations of even-even and even-odd $\pdo$s. First,
it holds in the following circumstances.
\begin{cor}\label{c:zetaglobaldefect}
\noi {\rm (i)} \ If $A$ satisfies one of the cases {\rm (1), (2)}
or {\rm (3)} of \propref{p:TRdefined} then ${\rm TR}(A)$ is
defined and \eqref{e:zetaglobaldefect} holds. In case (1) this
reduces to
\begin{equation}\label{e:zetatzerocases}
\zeta_\theta(A, Q, 0 ) \ =\  \TR(A)  \ - \ {\rm tr} (A\,\Pi_Q).
\end{equation}
If $Q$ is an even-even operator and has even order, then
\eqref{e:zetatzerocases} also holds when $A$ satisfies case {\rm
(2)} (assumes $M$ is odd-dimensional) or {\rm (3)} (assumes $M$ is
even-dimensional) of \propref{p:TRdefined}. These facts are known
from \cite{Gr1}.

\vskip 2mm \noi {\rm (ii)} \ If $A$ is a differential operator,
and more generally whenever $\TRx(A)=0$ for all $x\in M$,
\eqref{e:zetaglobaldefect} reduces to
\begin{equation}\label{e:TRAzero2}
\z_{\th}(A,Q,0) \ = \  - \ \frac{1}{q}\,{\rm res}(A\,\log_{\th}Q)
\ - \ {\rm tr} (A\,\Pi_Q).
\end{equation}
\end{cor}

\begin{proof}
Part (ii) follows from \propref{prop:cutoffdiff}. For part (i), it
is clear that \eqref{e:zetatzerocases} holds when ${\rm
res}_x(A\,\log_{\th}Q) =0$ for each $x\in M$. This is evident for
case {\rm (1)} operators. If $A$ satisfies case {\rm (2)} (resp.
case {\rm (3)}) of \propref{p:TRdefined} and if $Q$ is even-even
and of even order, then it is not hard to see that
$\s_{A\,\log_{\th}Q}(x,\xi)$ is also  even-even (resp. even-odd)
and hence $(\s_{A\,\log_{\th}Q})_{-n,0}(x,\xi)$ vanishes when
integrated over the $n-1$ sphere.
\end{proof}

\begin{ex}\label{rk:oddorder}
{\rm To see that \eqref{e:zetaglobaldefect} may hold with all
three terms non-zero, take $A = D + S$ with $D$ a differential
operator and $S$ a smoothing operator, and let $Q\in \Ell_{ord
>0}^{adm}(M, E)$. Then $\TR(A) = \tr(S)$ and $\res(A\,\log_{\th}Q)
= \res(D\,\log_{\th}Q)$ both exist (note
\corref{c:zetaglobaldefect} (ii)) and are non-zero in general. For
example, if $Q=D\in \Ell_{ord
>0}^{adm}(M, E)$ is invertible one has $\res(D\log_{\th}D) =
-\,\z_{\th}(D,-1)$.}
\end{ex}

\begin{rk}
In \corref{c:zetaglobaldefect} {\rm (i)}, if $Q$ has odd-order
then \eqref{e:zetaglobaldefect} may hold with all three terms
non-zero due to dependence on the choice of the spectral cut. The
distinct behaviour for odd-order $Q$ was kindly pointed out to the
authors by Gerd Grubb.
\end{rk}

\begin{rk} Using \thmref{prop:globaldefect}
similar facts to those in \corref{c:zetaglobaldefect} can be seen
to hold for the $\zeta^{(k)}_\theta(A, Q,\frac{\alpha-j}{q})$, see
also \cite{Gr1} Sect.3. The regularity of $\zeta_\theta(A, Q,z)$
at $z=0$ in (ii) is proved in \cite{GS}. When $A=I$ the identity
\eqref{e:TRAzero2} was shown in \cite{Sc}. On the other hand, when
$Q$ is a differential operator and taking $A=Q^m$ in
\eqref{e:TRAzero2} gives
\begin{equation}\label{eq:zetaatzerok}
\zeta_\theta( Q, -m) \ = \ -\, \frac{1}{q}\ {\rm res}\left(Q^m
\log_\theta Q\right) \ - \ {\rm tr} (Q^m\,\Pi_Q),
\end{equation}
which was obtained in the case when $Q$ is positive and invertible
by other methods in \cite{Lo}. Note that for sufficiently large
$m$ one has ${\rm tr} (Q^m\,\Pi_Q)=0$.
\end{rk}

\vskip 2mm

Looking at the next term up in the Laurent expansion, around $z=0$
the zeta function $\z_{\th}(Q,z) = \TR(Q_{\th}^{-z})$ is
holomorphic and hence the $\zeta$-determinant
$${\rm det}_{\zeta, \theta} Q \ = \exp(-\,\zeta_\theta^\prime(Q, 0)),$$
is defined, where $\zeta_\theta^\prime(Q, 0) =
\dd_z\zeta_\theta(Q, z))_{\vert_{z=0}}$.
\begin{thm}\label{t:det} One has
\begin{equation}\label{eq:logdet2}
\log {\rm det}_{\zeta, \theta} (Q) = \int_M dx\left( \TRx
\left(\log_\theta Q\right)-\frac{1}{ 2q}\, \resxo
\left(\log^2_\theta Q\right)\right) - \ {\rm tr} ( \log_\theta
Q\,\Pi_Q).
\end{equation}
If $M$ is odd-dimensional and $Q$ is an even-even operator and has
even order then one has (as known from \cite{O2},\cite{Gr1}
Sect.3, see also \cite{KV} Sect. 4)
\begin{equation}\label{eq:logdet2b}
\log {\rm det}_{\zeta, \theta} (Q) = \TR \left(\log_\theta
Q\right) - \ {\rm tr} ( \log_\theta Q\,\Pi_Q),
\end{equation}
where $\TR \left(\log_\theta Q\right)=\int_M  \TRx
\left(\log_\theta Q\right)\,dx,$
\end{thm}
\begin{proof}
Examining the coefficient of $z$ in the Laurent expansion
\eqref{eq:KV3zerotrace} immediately yields \eqref{eq:logdet2}. If
$Q$ is even-even and of even order then the classical component of
the local symbol of $\log^2_\theta Q\in\Cl^{0,2}(M,E)$ also has
even-even parity. Hence the local residue integral of the term of
homogeneity $-n$ then vanishes, $\TRx \left(\log_\theta
Q\right)\,dx$ defines a global density on $M$, and
\eqref{eq:logdet2} reduces to \eqref{eq:logdet2b}.
\end{proof}

\subsection{The canonical trace on commutators and the residue trace on logarithms}

The canonical trace $\TR$ is not defined on a commutator  of
classical $\pdo$s which has integer order. Rather the following
property holds.

\begin{thm}\label{prop:rescommlog}
Let $Q\in \Ell^{ adm}_{ord>0}(M, E)$ be of order $q$ and with
spectral cut $\theta$, and let $A\in \Cl^{\a}(M,  E)$, $B\in
\Cl^{\b}(M, E)$ for any $\a,\b\in\Rf$. Then
$$\left(\,\TRx\left([A,B]\right) - \frac{1}{q}\,{\rm
res}_{x,0}\left([A, B\log_\theta Q]\right)\,\right) dx$$ defines a
global density on $M$ and one has
\begin{equation}\label{e:TRacedefect}
\int_M dx\, \left(\,\TRx\left([A,B]\right) - \frac{1}{q}\,{\rm
res}_{x,0}\left([A, B\log_\theta Q]\right)\,\right) \ = \ 0
\end{equation}
independently of the choice of  $Q$.
\end{thm}
\begin{proof}
Using the vanishing of ${\rm TR}$ in \propref{p:TRtrace} (1), for
$z\neq 0$ sufficiently close to $0$ we have
\begin{equation}\label{e:TRABQz}
{\rm TR} \left( [A, B\,Q_\theta^{-z}]\right) = 0.
\end{equation}
Hence the function $z{\rm TR} \left( [A, B\,Q_\theta^{-z}]\right)$
also vanishes identically for such non-zero $z$. But from
\eqref{eq:KV3zerotrace}, $z{\rm TR} \left( [A,
B\,Q_\theta^{-z}]\right)$ extends holomorphically  to include
 $z=0$. By equation \eqref{e:TRABQz} this
analytically continued function must also vanish at $z=0$. It
follows that ${\rm TR} \left( [A, B\,Q_\theta^{-z}]\right)$ is
holomorphic near $z=0$ and so \eqref{eq:KV3zerotrace} implies
${\rm fp}_{z=0} {\rm TR} \left( [A,
B\,Q_\theta^{-z}]\right)=\lim_{z\to 0}{\rm TR} \left( [A,
B\,Q_\theta^{-z}]\right)=0.$ Applying Proposition
\ref{prop:globaldefect} to $A(z)= [A, B\,Q_\theta^{-z}]$ with
$z_0=0$  we have by \thmref{cor:defect}
\begin{eqnarray*}\label{eq:reslog}
0 & = & {\rm fp}_{z=0} {\rm TR} \left( [A,
B\,Q_\theta^{-z}]\right)\\ & =  & \int_M dx\, \left(\, \TRx \left(
[A, B (I- \Pi_Q)]\right)+ \frac{1}{q} \, \resxo\left( [A,
B\log_\theta Q]\right) \,\right)
\end{eqnarray*}
which is equation \eqref{e:TRacedefect}, since ${\rm TR} \left(
[A, B\Pi_Q)]\right) = {\rm tr} \left( [A, B\Pi_Q)]\right) = 0.$
\end{proof}

\begin{cor}
Let $Q\in \Ell^{ adm}_{ord>0}(M, E)$ be of order $q$ and with
spectral cut $\theta$, and let $A\in \Cl^{\a}(M,  E)$, $B\in
\Cl^{\b}(M, E)$. Then in cases {\rm (1), (2)} and {\rm (3)} of
\propref{p:TRtrace} the form ${\rm res}_x\left([A, B\log_\theta
Q]\right)\,dx$ determines a global density on $M$ and one has
\begin{equation*}
{\rm res}\left([A, B\log_\theta Q]\right) =0
\end{equation*}
independently of the choice of  $Q$.
\end{cor}

\begin{rk}
The independence from $Q$ can also be seen for the residue trace
term directly; given $Q_1, Q_2\in \Ell^{ adm}_{ord>0}(M, E)$ of
order $q_1$ and $q_2$ respectively with common spectral cut
$\theta$, the difference
$$\left(\frac{1}{q_1} {\rm res}_{x,0}\left( [A, B\log_\theta
Q_1]\right)- \frac{1}{q_2} {\rm res}_{x,0}\left( [A, B\log_\theta
Q_2]\right)\right)\, dx$$ defines a global density  which
integrates to
$$ {\rm res}\left( \left[A, B \left(\frac{\log_\theta Q_1}{q_1}-
\frac{\log_\theta Q_2}{q_2}\right)\right]\right)=  0$$ since
$\frac{\log_\theta Q_1}{q_1}- \frac{\log_\theta Q_2}{q_2}$ is a
classical $\pdo$.
\end{rk}

\vskip 2mm

A useful  consequence of \thmref{prop:rescommlog} and
\propref{prop:cutoffdiff} is:
\begin{cor}
\noi Let $Q\in \Ell^{ adm}_{ord>0}(M, E)$ of order $q$ and with
spectral cut $\theta$ and  let $A, B\in \Cl(M,  E)$. Whenever
${\rm TR}\left([A, B]\right)=\int_M\, dx\,{\rm TR}_x\left([A,
B]\right)  $ is  well defined  then ${\rm res}_{x,0}\left([A,
B\log_\theta Q]\right) dx$ is globally defined and one then has
\begin{equation}\label{eq:TRbracketres}
{\rm res}\left( [A, B\log_\theta Q]\right)=q\, {\rm
TR}\left([A,B]\right).
\end{equation}
In particular, if  $[A, B]$ is a differential operator then ${\rm
res}_{x,0}\left([A, B\log_\theta Q]\right) dx$ is globally defined
and one has
$${\rm res}\left( [A, B\log_\theta Q]\right)=0.$$

\vskip 1mm

\noi In that case, whenever \ ${\rm res}_{x,0}\left(A\,
B\,\log_\theta Q\right) dx$ defines a global density, then so does
${\rm res}_{x,0}\left(B\, \log_\theta Q\, A\right) dx$ and
$$  {\rm res}\left( B\, \log_\theta Q\, A\right)=
{\rm res}\left( A\,  B\, \log_\theta Q\right). $$ In particular,
since  ${\rm res}\left(  \log_\theta Q\right)$ exists  \cite{O1},
for any invertible $A\in \Cl(M,  E)$
\begin{equation}\label{e:logtracedefect}
{\rm res}\left( A^{-1}\, \log_\theta Q\, A\right)=  {\rm
res}\left( \log_\theta Q\right).
\end{equation}
\end{cor}

\begin{rk}
This proposition partially generalizes the fact \cite{O1} that  $
{\rm res}_{x,0}\left([A, \log_\theta Q]\right) dx$ for $A$ a
classical $\pdo$ defines a global density and  ${\rm res}\left([A,
\log_\theta Q]\right) =0$, which when $A$ is a differential
operator follows from the corollary applied to $B=I$.
\end{rk}

On the other hand, the well known (\cite{MN}, \cite{O1},
\cite{CDMP}, \cite{Gr2}) trace defect formula
\begin{equation}\label{e:TRacedefect2a}
\z_{\th}([A, B],Q,0)  =  - \frac{1}{q} \,{\rm res}\left( A\,[ B,
\log_\theta Q]\right).
\end{equation}
follows easily by applying the same argument as in the proof of
\thmref{prop:rescommlog} to $C(z)= A\, [B, Q^{-z}]$. From
\eqref{e:TRacedefect} and \eqref{e:TRacedefect2a} we infer:

\begin{cor}
For classical $\pdo$s $A$ and $B$
\begin{equation*}
- \frac{1}{q} \,{\rm res}\left( A\,[ B, \log_\theta Q]\right) \ =
\ \int_M dx\, \left(\TRx([A, B]) \ -\  \frac{1}{q} \,{\rm
res}_x\left( [A , B]\log_\theta Q \right)\,\right).
\end{equation*}
In cases {\rm (1), (2)} and {\rm (3)} of \propref{p:TRtrace} the
form ${\rm res}_x\left( [A , B]\log_\theta Q \right)\,dx$
determines a global density on $M$ and one has
\begin{equation*} {\rm res}\left( A\,[ B, \log_\theta Q]\right) \
= \,{\rm res}\left( [A , B]\log_\theta Q \right).
\end{equation*}

\end{cor}

\vskip 3mm

\noi While from \propref{prop:rescommlog} we conclude:

\begin{cor}\label{cor:reslogbrackets}
The density $\resx\left( [A , B\log_\theta Q]  -  [A ,
B]\log_\theta Q \right)\,dx$ is globally defined on $M$ for
classical $\pdo$s $A$ and $B$ and one has
\begin{equation*}
{\rm res}\left( A\,[ B, \log_\theta Q]\right) \ = \ \res\left( [A
, B]\log_\theta Q- [A , B\log_\theta Q] \right)
\end{equation*}
\end{cor}
\begin{proof}
\begin{eqnarray*}
\frac{1}{q} \,{\rm res}\left( A\,[ B, \log_\theta Q]\right) & =& -
\int_M dx\, \left(\TRx([A, B]) \ -\  \frac{1}{q} \,{\rm res}_x( [A
, B]\log_\theta Q) \right)\\
&=& - \int_M dx\, \left(\TRx([A, B]) \ -\  \frac{1}{q} \,{\rm
res}_x\left( [A , B\, \log_\theta Q ]\right)\right.\\
&& + \left. \frac{1}{q} \,{\rm res}_x\left(  [A , B]\log_\theta
Q-[A , B\, \log_\theta Q ]\right)\right)\\
&=&   \frac{1}{q} \int_M dx\,\,{\rm res}_x\left(  [A ,
B]\log_\theta Q-[A , B\, \log_\theta Q ]\right)\\
&=&\frac{1}{q} \ {\rm res}\left(  [A , B]\log_\theta Q-[A , B\,
\log_\theta Q ]\right)\\
\end{eqnarray*}
\end{proof}

\vskip 3mm

We point out that \corref{cor:index1} and \eqref{e:TRacedefect2a}
imply the following local index formulae.
\begin{cor}\label{cor:index1}
Let $A$ be an elliptic $\pdo$ with parametrix $B$. Let $Q\in
\Ell^{ adm}_{ord>0}(M, E)$ be of order $q$ and with spectral cut
$\theta$. Then, independently of the choice of  $Q$,

\begin{equation}\label{e:resbracketlog}
{\rm res}\left( [A, B\log_\theta Q]\right) \ = \ {\rm res}\left(
A\,[ B, \log_\theta Q]\right)
\end{equation}
and are equal to $-q\,\index(A)$.
\end{cor}
\begin{proof}
In this case $\index(A) = \tr([A,B])$ and since $[A, B]$ is
smoothing equal to ${\rm TR}\left([A, B]\right)$. The first
equality thus follows from \eqref{eq:TRbracketres}. Since $AB = I
+ S$ where $S$ is a smoothing operator, and since
$\resxo(S\log_\theta Q)$ is therefore equal to zero, the second
equality also follows.
\end{proof}

\section*{Appendix A: Proof of the density formula}

The purpose here is to give a direct elementary proof of
\thmref{t:gendensity}, which for the family $z\mapsto A(z)\in
\Cl(M, E)$ parametrized by a domain $W\subset \C$ states that
irrespective of the order $\a(z_0)\in\Rf$ of $A(z_0)$
\begin{equation}\label{e:globaldensity}
\left(\TRx(A) \ - \ \frac{1}{\alpha^\prime}\,
\resxo(A^\prime)\,\right)\, dx
\end{equation}
defines a global density on $M$. Here, we have written $A =
A(z_0)$,   $A^\prime = A^\prime(z_0) := d/dz|_{z=z_0}(A(z))$,
and $\a^\prime = \a^\prime(z_0)$.

\vskip 2mm

\noi From previous works \cite{KV} it is known that $\TRx(A(z_0))
\, dx$ defines a global density on $M$ when $\a(z_0)$ is not
integer valued; this follows immediately from
\eqref{e:globaldensity} and \propref{prop:diffsymb}.

\vskip 3mm

\noi The method of proof uses a generalization of the method used
in \cite{O1} to show that the residue density is globally defined
for any classical $\pdo$, and the method in \cite{Le} used to show
that the canonical density is globally defined for classical
$\pdo$s of non-integer order. We will take $A$ to be scalar valued
for notational brevity, but the proof works in the same way for
endomorphism valued operators; indeed it works equally for the
pre-tracial density $ (\cutoffint_{T_xM}\s_{A}(x,\xi)\,\dbar\xi \
- \ \frac{1}{\alpha^\prime}\,
\int_{S^*_xM}\,(\s_{A^\prime})_{-n,0} (x, \xi) \,\dbar_S\xi )\,
dx. $

\vskip 3mm

First, we have a lemma, generalizing Lemma C.1 in \cite{O1}.

\begin{lem}
Let $f(\xi)$ be a smooth function on $\Rf^n$ which is homogeneous
of degree $-n$ for $|\xi|\geq 1$ and let $T$ be an invertible
linear map on $\Rf^n$. Then for $s\in\Cf$ and any non-negative
integer $k$

\begin{equation*}\label{e:homogonsphere}
\int_{|\eta|=1} f(T\eta) \,|T\eta|^s\, \log^k|T\eta| \,
\dbar_S\eta \ = \ \frac{(-1)^k}{|{\rm det} T|} \,\int_{|\xi|=1}
f(\xi)\,|T\ii\xi|^{-s}\, \log^k|T\ii\xi| \,\dbar_S\xi.
\end{equation*}

\vskip 1mm

\noi Specifically, one has

\begin{equation}\label{e:homogonsphere2}
\int_{|\eta|=1} f(T\eta) \, \log|T\eta| \, \dbar_S\eta \ = \
\frac{-1}{|{\rm det} T|} \,\int_{|\xi|=1} f(\xi) \, \log|T\ii\xi|
\,\dbar_S\xi.
\end{equation}

\begin{equation}\label{e:homogonsphere3}
\int_{|\eta|=1} f(T\eta) \, \dbar_S\eta \ = \ \frac{1}{|{\rm det}
T|} \,\int_{|\xi|=1} f(\xi) \,\dbar_S\xi.
\end{equation}

\vskip 1mm

\end{lem}

\begin{proof}
It is enough to prove this for $k=0$, differentiation with respect
to $s$ yields the general formula. We have, using the linearity of
$T$,
\begin{eqnarray}
\int_{1\leq|\eta|\leq 2} f(T\eta)\,|T\eta|^s \, d\eta & = &
\int_{|\eta|=1}\int_{1\leq r \leq 2} f(rT\eta)\,r^s |T\eta|^s \,
r^{n-1}\,dr\,\dbar_S\eta \nonumber\\[3mm]
& = & \left(\frac{2^s - 1}{s}\right)\,\int_{|\eta|=1} f(T\eta)\,
|T\eta|^s \,
\dbar_S\eta.\nonumber % \label{e:firsteval}
\end{eqnarray}

\noi On the other hand, changing variable,

\begin{eqnarray*}
\int_{1\leq|\eta|\leq 2} f(T\eta)\,|T\eta|^s \, d\eta & = &
\frac{1}{|{\rm det} T|} \,\int_{1\leq |T\ii\eta|\leq 2}
f(\eta)\,|\eta|^s \, d\eta \nonumber\\[3mm]
& = & \frac{1}{|{\rm det} T|}
\,\int_{|\eta|=1}\int_{1/|T\ii\eta|\leq r \leq  2/|T\ii\eta|}
f(r\eta)\,r^s |\eta|^s \,
r^{n-1}\,dr\,\dbar_S\eta \nonumber\\[3mm]
& = & \frac{1}{|{\rm det} T|} \, \left(\frac{2^s - 1}{s}\right)\,
\int_{|\eta|=1} f(\eta)\, |T\ii\eta|^{-s}\,\dbar_S\eta.
\nonumber%\label{e:secondeval}
\end{eqnarray*}
\end{proof}

\vskip 2mm

Consider now a local chart on $M$ defined by a diffeomorphism $x :
\Omega \too U$ from an open subset $\Omega$ of $M$ to an open
subset $U$ of $\Rf^n$. For $p\in\Omega$ we then have the local
coordinate $x(p)\in\Rf^n$. Let $\k : U \too V$ be a diffeomorphism
to a second open subset $V$ of $\Rf^n$. Then $y(p) = \k(x(p))$ is
also a local coordinate for $\Omega$.

\vskip 2mm

\noi Let $a(x(p),\,\xi) = \tilde{a}(x(p),x(p),\xi) $ where
$\tilde{a}(x(p),y(p),\xi)$ denotes the local amplitude of $A$ in
$x$-coordinates, and likewise let $b(y(p),\xi)$ denote the
amplitude along the diagonal in $y$-coordinates. From \cite{Ho}
with $T(p):= (D\k_{x(p)})^t$ we have
\begin{eqnarray*}
{\rm TR}_{y(p)}(A)\,dy(p) & := & \cutoffint_{\Rf^n}
b(y(p),\,\xi)\,\dbar\xi\,dy(p) \nonumber \\ & = &
\cutoffint_{\Rf^n}a(x(p),\,T(p)\xi)\,\dbar\xi\,dy(p).
\label{e:globaldensity3}
\end{eqnarray*}
According to the transformation rule in \propref{prop:noninv}, for
$f\in \CS(V)$ and $T$ an invertible linear map on $\Rf^n$
\begin{equation*}\label{e:classicaltransform}
\cutoffint_{\Rf^n} f(T\xi) \, \dbar\xi = \frac{1}{|{\rm det}
T|}\left(\cutoffint_{\Rf^n} f(\xi) \, \dbar\xi - \int_{|\xi|=1}
f(\xi)_{(-n)}\,\log|T\ii \xi| \, \dbar\xi\right)
\end{equation*}
with $f(\xi)_{(-n)}$ the homogeneous component of $f$ of degree
$-n$. Hence
\begin{eqnarray}
& & {\rm TR}_{y(p)}(A) \,dy(p) \nonumber\\ & & = \ \
\frac{1}{|{\rm det} T(p)|}\left(\cutoffint_{\Rf^n}
a(x(p),\xi)\,\dbar\xi\,dy(p) - \int_{|\xi|=1}
a(x(p),\xi)_{(-n)}\,\log|T(p)\ii \xi| \,
\dbar\xi \,dy(p)\right) \nonumber\\
 & & = \ \ \cutoffint_{\Rf^n} a(x(p),\xi)\,\dbar\xi\,dx(p)  -
\int_{|\xi|=1} a(x(p),\xi)_{(-n)}\,\log|T(p)\ii \xi| \,
\dbar\xi\,dx(p)\nonumber\\
 & & = \ \ {\rm TR}_{x(p)}(A) \,dx(p)  -
\int_{|\xi|=1} a(x(p),\xi)_{(-n)}\,\log|T(p)\ii \xi| \,
\dbar\xi\,dx(p).\label{e:TRtransform}
\end{eqnarray}

\vskip 2mm

We turn now to the other component of \eqref{e:globaldensity}
given in $y$-coordinates by
\begin{equation*}\label{e:localresdensity}
 - \ \frac{1}{\alpha^\prime}\, \int_{|\xi|=1}
 b^\prime(y(p),\xi)_{(-n)}\,\dbar_S\xi\,
 dy(p),
\end{equation*}
where $b^\prime(y(p),\xi) = d/dz|_{z=z_0}(\s_A(z)(y(p),\xi))$ is
the symbol derivative in $y$-coordinates and where
$b^\prime(y(p),\xi)_{(-n)}$ denotes its log-homogeneous (cf.
\eqref{eq:logpolysymb}) component of degree $-n$. From \cite{Ho}
we have the asymptotic formula

\begin{equation}\label{e:symboltransform}
b^\prime(y(p),\xi) \sim \sum_{|\mu|\geq
0}\dd^{\mu}_{\xi}a^\prime(x(p),T(p)\xi)\Psi_{\mu}(x,\xi)
\end{equation}
with $\Psi_{\mu}(x,\xi)$ polynomial in $\xi$ of degree of at most
$|\a|/2$. To begin with, suppose that $a(x(p),\xi)$ is homogeneous
in $\xi$ of degree $-n$. Then from \eqref{e:derivsigma} for
$|\eta|\geq 1$

\begin{equation}\label{e:aprimehomog}
a^\prime(x(p), \eta) = \alpha^\prime\,
a(x(p), \eta)\,\log |\eta| +  p_{-n}(x(p), \eta)\\[1mm]
\end{equation}

\noi with $p_{-n}(x(p), \eta)$ positively homogeneous in $\eta$ of
degree $-n$, and $a^\prime(x(p), \eta) = a^\prime(x(p),
\eta)_{(-n)}$. Thus, if $a(x(p),\xi)$ is homogeneous in $\xi$ of
degree $-n$, by \eqref{e:symboltransform} and
\eqref{e:aprimehomog}

\begin{eqnarray}
- \ \frac{1}{\alpha^\prime}\, \int_{|\xi|=1}
 b^\prime(y(p),\xi)_{(-n)}\,\dbar_S\xi\,
 dy(p) &  = & - \ \frac{1}{\alpha^\prime}\, \int_{|\xi|=1}
 a^\prime(x(p),T(p)\xi)\,\dbar_S\xi\,
 dy(p) \nonumber\\
& = & - \ \frac{1}{\alpha^\prime}\, \int_{|\xi|=1}
 \alpha^\prime\,
a(x(p), T(p)\xi)\,\log |T(p)\xi| \,\dbar_S\xi\,
 dy(p)   \nonumber\\ & & \hskip 10mm -
\frac{1}{\alpha^\prime}\, \int_{|\xi|=1} p_{-n}(x(p),
T(p)\xi)\,\dbar_S\xi\,
 dy(p)  \nonumber\\
 & = & - \ \int_{|\xi|=1}
a(x(p), T(p)\xi)\,\log |T(p)\xi| \,\dbar_S\xi\,
 dy(p)   \nonumber\\ & & \hskip 10mm -
\frac{1}{\alpha^\prime}\, \int_{|\xi|=1} p_{-n}(x(p),
T(p)\xi)\,\dbar_S\xi\,
 dy(p). \label{e:localresdensity1}
\end{eqnarray}
Using equations \eqref{e:homogonsphere2} and
\eqref{e:homogonsphere3} of \lemref{e:homogonsphere},
\eqref{e:localresdensity1} becomes
%$a^\p$ so that \eqref{e:symboltransform} gives $
%b^\prime(y(p),\xi)_{(-n)} = a(x(p),\xi)$. Then we have

\begin{eqnarray}
- \ \frac{1}{\alpha^\prime}\, {\rm \res}_{y(p),0}(A^\prime)
 dy(p) & & = \ \frac{1}{|{\rm det} T(p)|} \,  \int_{|\xi|=1} a(x(p),
\xi)\,\log |T(p)\ii\xi| \,\dbar_S\xi\,
 dy(p) \nonumber\\ & &    \hskip 30mm -
\frac{1}{\alpha^\prime}\, \frac{1}{|{\rm det} T(p)|} \,
\int_{|\xi|=1} p_{-n}(x(p),\xi)\,\dbar_S\xi\,
 dy(p)\nonumber\\[3mm]
& & =    \ \int_{|\xi|=1} a(x(p), \xi)\,\log |T(p)\ii\xi|
\,\dbar_S\xi\,
 dx(p) \nonumber\\ & &    \hskip 30mm -
\frac{1}{\alpha^\prime}\, \int_{|\xi|=1}
p_{-n}(x(p),\xi)\,\dbar_S\xi\,
 dx(p) \nonumber\\[3mm]
& & =  \ \int_{|\xi|=1} a(x(p), \xi)\,\log |T(p)\ii\xi|
\,\dbar_S\xi\,
 dx(p) \nonumber\\
  & &    \hskip 30mm - \ \frac{1}{\alpha^\prime}\, {\rm \res}_{x(p),0}(A^\prime)
 dx(p),\label{e:localresdensity2}
\end{eqnarray}

\noi where the final equality follows from \eqref{e:aprimehomog}.
Adding \eqref{e:TRtransform} and \eqref{e:localresdensity2} we
have when $a(x(p),\xi)$ is homogeneous in $\xi$ of degree $-n$

\begin{equation}\label{e:globaldensity2}
\left({\rm TR}_{y(p)}(A) -  \frac{1}{\alpha^\prime(z_0)}\, {\rm
\res}_{y(p),0}(A^\prime)\right) dy(p) = \left({\rm TR}_{x(p)}(A) -
\frac{1}{\alpha^\prime(z_0)}\, {\rm
\res}_{x(p),0}(A^\prime)\right) dx(p),
\end{equation}

\noi proving the invariance of \eqref{e:globaldensity} in this
case.

Next suppose that $a(x(p),\xi)$ is homogeneous in $\xi$ of degree
$\a > -n$. Then from \eqref{e:symboltransform} and since we can
commute the $z$ and $\mu$ derivatives
\begin{equation*}
b^\prime(y(p),\xi)_{(-n)} = \sum_{|\mu|\geq \a +
n}\left.\frac{d}{dz}\right|_{z=z_0}\,\dd^{\mu}_{\xi}
\left(a(z)(x(p),T(p)\xi)\right)\,\Psi_{\mu,-n}(x,\xi).
\end{equation*}
where $\Psi_{\mu,-n}(x,\xi)$ is a polynomial in $\xi$ of degree
$|\mu|-n-\a$. Hence
\begin{eqnarray*}
& & \int_{|\xi|=1}
 b^\prime(y(p),\xi)_{(-n)}\,\dbar_S\xi\,
 dy(p) \\& & \hskip 10mm =  \sum_{|\mu|\geq \a +
n}\left.\frac{d}{dz}\right|_{z=z_0}\, \int_{|\xi|=1}
 \dd^{\mu}_{\xi}\left(a(z)(x(p),T(p)\xi)\right)\,\Psi_{\mu,-n}(x,\xi)\,\dbar_S\xi\,
 dy(p) \nonumber\\
& & \hskip 10mm  = 0.
\end{eqnarray*}
The final equality follows using the integration by parts property
in Lemma C1 of \cite{O1}, which states that if $g(\xi)$ and
$h(\xi)$ are homogeneous in $\xi$ of degrees $\g, \delta$ where
$\g + \delta = 1-n$, then
$$ \int_{|\xi|=1}
 (\dd_{\xi_j}g(\xi))h(\xi)\,\dbar_S\xi  = - \int_{|\xi|=1}
 g(\xi)\dd_{\xi_j}h(\xi)\,\dbar_S\xi,$$
along with the fact that $\Psi_{\mu,-n}(x,\xi)$ polynomial in
$\xi$ of degree $|\mu|-n-\a$.

\vskip 2mm

This completes the proof that \eqref{e:globaldensity} is a density
independent of coordinates.

\section*{Appendix B: Proof of \lemref{lem:cutoff} and \lemref{cutoffint}}
For a fixed $N\in \N$ chosen large enough such that
$\re(\alpha)-N-1<-n$, we write $\sigma(x,\xi)= \sum_{j=0}^{K_N}
\sigma_{\alpha-j} (x,\xi)+\sigma_{(N)}(x,\xi)$ and split the
integral accordingly as
$$\int_{B_x^*(0, R)}  \sigma(x,\xi) \dbar\xi=
\sum_{j=0}^{N}\int_{B_x^*(0, R)} \sigma_{\alpha-j}(x,\xi)
\dbar\xi+ \int_{B_x^*(0, R)} \sigma_{(N)}(x,\xi) \dbar\xi.$$ Since
$\re(\alpha)-N-1<-n$,  $\sigma_{(N)}$ lies in $ L^1(T_x^*U)$ and
the integral  $\int_{B_x^*(0, R)} \sigma_{(N)}(x,\xi) \dbar\xi$
converges when $R\to \infty$ to $\int_{T_x^*U} \sigma_{(N)}(x,\xi)
\dbar\xi$. On the other hand, for any $j\leq N$
\begin{eqnarray}\label{eq:intlambda}
\int_{B_x^*(0, R)} \sigma_{\alpha-j}= \int_{B_x^*(0, 1)}
\sigma_{\alpha-j}+ \int_{D_x^*(1,R)} \sigma_{\alpha-j}.
\end{eqnarray}
Here $D_x^*(1, R)= B_x^*(0,R)\backslash B_x^*(0, 1)$. The first
integral on the r.h.s. converges and since
$\sigma_{\alpha-j}(x,\xi)\sim\sum_{l=0}^k \sigma_{\alpha-j,
l}(x,\xi)\log^l [\xi],$  the second integral reads:
$$\int_{D_x^*(1,R)} \sigma_{\alpha-j}(x,\xi) \dbar\xi=
\sum_{l=0}^k\int_1^R  r^{\alpha-j+n-1}\log^l r \, dr\cdot
\int_{S_x^*U} \sigma_{\alpha-j,l}(x, \omega)d\omega.$$ Hence the
following asymptotic behaviour:
\begin{eqnarray*}
&{}& \int_{D_x^*(1,R)} \dbar\xi\, \sigma_{\alpha-j}(x,\xi) \ \ \sim_{R\to \infty}\\
& & \sum_{l=0}^k\frac{ \log^{l+1}R}{l+1}\cdot \int_{S_x^*U}
\sigma_{\alpha-j,l}( x,\omega)\dbar_S\xi=\sum_{l=0}^k
\frac{\log^{l+1}R}{l+1} \int_{S_x^*U} \sigma_{-n,l } (x,\xi) \
\,\dbar_S\xi\quad \, {\rm if} \, \, \alpha-j=-n
\end{eqnarray*}
\begin{eqnarray*}
\int_{D_x^*(1,R)}\dbar\xi\, \sigma_{\alpha-j}(x,\xi)
 &\sim_{R\to \infty}& \sum_{l=0}^k
\left(\sum_{i=0}^l\frac{(-1)^{i+1}\frac{l!}{(l-i)!}\log^i
R}{(\alpha-j+n)^i}\cdot R^{\alpha-j+n}\int_{S_x^*U}
\sigma_{\alpha-j,l}(x, \xi)\,\dbar_S\xi\right.\\ &+&
(-1)^{l}l!\frac{R^{\alpha-j+n}}{(\alpha-j+n)^{l+1}}\cdot
 \int_{S_x^*U} \sigma_{\alpha-j,l} (x,\xi) \,\dbar_S\xi \\
&+& \left. \frac{(-1)^{l+1}l!}{(\alpha-j+n)^{l+1}}\cdot
\int_{S_x^*U} \sigma_{\alpha-j,l} (x,\xi) \,\dbar_S\xi\right)\,
{\rm
if}\quad \, \, \alpha-j\neq -n.\\
\end{eqnarray*}
Putting together these asymptotic expansions  yield the statements
with
$$C_x (\sigma)=
\int_{T_x^*U}\sigma_{(N)}+\sum_{j=0}^{N}  \int_{B_x^*(0, 1)}
\sigma_{a_j}+ \sum_{j=0, a_j+n\neq 0}^{N}\sum_{l=0}^L
\frac{(-1)^{l+1}l!}{(a_j+n)^{l+1}}\int_{S_x^*U} \sigma_{a_j,l}.$$
The $\mu$-dependence follows from
\begin{eqnarray*}
\log^{l+1}(\mu \, R)&=&\log^{l+1}R\, \left(1+\frac{\log\mu}{\log
\, R}\right)^{l+1}\\
&\sim_{R\to\infty}&\log^{l+1}R\,\sum_{k=0}^{l+1} C_{l+1}^k
\left(\frac{\log\mu}{\log R}\right)^k.\\
\end{eqnarray*}
The logarithmic terms $\sum_{l=0}^k \frac{1}{l+1}\int_{S_x^*U}
\sigma_{-n,l } (x,\xi) \ \dbar_S\xi\log^{l+1}(\mu\, R)$  therefore
contribute to the finite part by $\sum_{l=0}^k
\frac{\log^{l+1}\mu}{l+1}\cdot \int_{S_x^*U} \sigma_{-n,l }
(x,\xi) \ \dbar_S\xi$ as claimed in the lemma.

\bibliographystyle{plain}

\vskip 10mm {\Small \noi\textsc{Laboratoire de Math\'ematiques,
Complexe des C\'ezeaux, Universit\'e Blaise Pascal, 63 177
Aubi\`ere Cedex F. E-mail:} sylvie.paycha@math.univ-bpclermont.fr

\vskip 5mm

\noi \textsc{Department of Mathematics, King's College London.
 E-mail:} sgs@mth.kcl.ac.uk}

\end{document}